\documentclass[11pt]{article}
\pdfoutput=1

\usepackage[dvipsnames]{xcolor}

\def\beq{\begin{equation} }\def\eeq{\end{equation} }\def\1{\mathbf{1}}

\usepackage[framemethod=default]{mdframed}
\usepackage{caption}

\usepackage{indentfirst}
\usepackage{bm, mathrsfs, graphics,float,amssymb,amsmath,subeqnarray,setspace,graphicx,amsthm,epstopdf,subfigure, enumerate, color}
\usepackage[utf8]{inputenc}
\usepackage[colorlinks,
            linkcolor=red,
            anchorcolor=blue,
            citecolor=blue
            ]{hyperref}
\usepackage{natbib}

\usepackage{fullpage}
\numberwithin{equation}{section}

\newtheorem{lemma}{Lemma}
\newtheorem{theorem}{Theorem}

\newtheorem{definition}{Definition}
\newtheorem{corollary}[theorem]{Corollary}

\ifx\assumption\undefined
\newtheorem{assumption}{Assumption}
\fi

\newcommand{\x}{\mathbf{x}}
\newcommand{\y}{\mathbf{y}}

\renewcommand{\u}{\mathbf{u}}
\renewcommand{\v}{\mathbf{v}}

\newcommand{\X}{\mathbf{X}}
\newcommand{\Y}{\mathbf{Y}}
\newcommand{\Z}{\mathbf{Z}}
\newcommand{\A}{\mathbf{A}}
\newcommand{\B}{\mathbf{B}}

\newcommand{\D}{\mathbf{D}}

\newcommand{\I}{\mathbf{I}}
\newcommand{\U}{\mathbf{U}}

\newcommand{\V}{\mathbf{V}}
\renewcommand{\H}{\mathbf{H}}

\newcommand{\Q}{\mathbf{Q}}
\newcommand{\R}{\mathbf{R}}

\newcommand{\diag}{\mbox{diag}}

\newcommand{\rank}{\mbox{rank}}

\newcommand{\argmin}{\mbox{argmin}}

\newcommand{\tr}{\mbox{trace}}
\newcommand{\<}{\left\langle}
\renewcommand{\>}{\right\rangle}

\usepackage{bbm}

\usepackage{multirow}
\usepackage{tablefootnote}
\usepackage{colortbl}
\usepackage{hhline}

\usepackage{algorithm}
\usepackage{algorithmic}

\usepackage{mathtools}

\begin{document}
\title{
Provable Accelerated Gradient Method for Nonconvex Low Rank Optimization
}

\author{
Huan Li
\thanks{
Peking University;
email: lihuanss@pku.edu.cn; zlin@pku.edu.cn
}
\qquad
Zhouchen Lin
\footnotemark[1]
}

\maketitle

\begin{abstract}
Optimization over low rank matrices has broad applications in machine learning. For large scale problems, an attractive heuristic is to factorize the low rank matrix to a product of two much smaller matrices. In this paper, we study the nonconvex problem $\min_{\U\in\mathcal{R}^{n\times r}} g(\U)=f(\U\U^T)$ under the assumptions that $f(\X)$ is restricted $\mu$-strongly convex and $L$-smooth on the set $\{\X:\X\succeq 0,\rank(\X)\leq r\}$. We propose an accelerated gradient method with alternating constraint that operates directly on the $\U$ factors and show that the method has local linear convergence rate with the optimal dependence on the condition number of $\sqrt{L/\mu}$. Globally, our method converges to the critical point with zero gradient from any initializer. Our method also applies to the problem with the asymmetric factorization of $\X=\widetilde\U\widetilde\V^T$ and the same convergence result can be obtained. Extensive experimental results verify the advantage of our method.
\end{abstract}

\section{Introduction}

Low rank matrix estimation has broad applications in machine learning, computer vision and signal processing. In this paper, we consider the problem of the form:
\begin{eqnarray}
\min_{\X\in\mathcal{R}^{n\times n}} f(\X),\quad s.t.\quad \X\succeq 0,\label{convex_problem}
\end{eqnarray}
where there exists minimizer $\X^*$ of rank-$r$. We consider the case of $r\ll n$. Optimizing problem (\ref{convex_problem}) in the $\X$ space often requires computing at least the top-$r$ singular value/vectors in each iteration and $O(n^2)$ memory to store a large $n$ by $n$ matrix, which restricts the applications with huge size matrices. To reduce the computational cost as well as the storage space, many literatures exploit the observation that a positive semidefinite low rank matrix can be factorized as a product of two much smaller matrices, i.e., $\X=\U\U^T$, and study the following nonconvex problem instead:
\begin{eqnarray}
\min_{\U\in\mathcal{R}^{n\times r}} g(\U)=f(\U\U^T).\label{non_convex_problem0}
\end{eqnarray}
A wide family of problems can be cast as problem (\ref{non_convex_problem0}), including matrix sensing \citep{Bhojanapalli-2016}, matrix completion \citep{jain-2013}, one bit matrix completion \citep{Davenport-2014}, sparse principle component analysis \citep{Cai-2013} and factorization machine \citep{lin-2016}. In this paper, we study problem (\ref{non_convex_problem0}) and aim to propose an accelerated gradient method that operates on the $\U$ factors directly. The factorization in problem (\ref{non_convex_problem0}) makes $g(\U)$ nonconvex, even if $f(\X)$ is convex. Thus, proving the acceleration becomes a harder task than the analysis for convex programming.

\subsection{Related Work}

Recently, there is a trend to study the nonconvex problem (\ref{non_convex_problem0}) in the machine learning and optimization community. Recent developments come from two aspects: (1). The geometric aspect which proves that there is no spurious local minimum for some special cases of problem (\ref{non_convex_problem0}), e.g., matrix sensing \citep{Bhojanapalli-2016}, matrix completion \citep{Ge-2016}, and \cite{ge-2017,li-2018,zhu-2018} for a unified analysis. (2). The algorithmic aspect which analyzes the local linear convergence of some efficient schemes such as the gradient descent method. Examples include \citep{Burer2003,Burer2005,Boumal-2016,Tu-2016,zhangQ-2015,park-2016} for semidefinite programs, \citep{sun-2015,park-2013,hardt-2014,zheng-2016,zhao-2015} for matrix completion, \citep{zhao-2015,park-2013} for matrix sensing and \citep{Chen-2016} for Robust PCA. The local linear convergence rate of the gradient descent method is proved for problem (\ref{non_convex_problem0}) in a unified framework in \cite{Srinadh-2016-colt,Chen-2015,wang-2016}. However, no acceleration scheme is studied in these literatures. It remains an open problem on how to analyze the accelerated gradient method for nonconvex problem (\ref{non_convex_problem0}).

Nesterov's acceleration technique \cite{Nesterov1983,Nesterov1988,Nesterov-2004} has been empirically verified efficient on some nonconvex problems, e.g., Deep Learning \citep{Sutskever-2013}. Several literatures studied the accelerated gradient method and the inertial gradient descent method for the general nonconvex programming \cite{Saecd2013,li-2016-nips,Xu-2014}. However, they only proved the convergence and had no guarantee on the acceleration for nonconvex problems.  Carmon et al. \cite{Carmon-2016,Carmon-2017}, Agarwal et al. \cite{Agarwal-2017} and Jin et al. \cite{jin-2018} analyzed the accelerated gradient method for the general nonconvex optimization and proved the complexity of $O(\epsilon^{-7/4}\mbox{log}(1/\epsilon))$ to escape saddle points or achieve critical points. They studied the general problem and did not exploit the specification of problem (\ref{non_convex_problem0}). Thus, their complexity is sublinear.  Necoara et al. \cite{Necoara-2016} studied several conditions under which the gradient descent and accelerated gradient method converge linearly for non-strongly convex optimization. Their conclusion of the gradient descent method can be extended to nonconvex problem (\ref{non_convex_problem0}). For the accelerated gradient method,  Necoara et al. required a strong assumption that all $\y^k,k=0,1,\cdots,$\footnote{Necoara et al. \cite{Necoara-2016} analyzed the method with recursions of  $\y^k=\x^k+\frac{\sqrt{L}-\sqrt{\mu}}{\sqrt{L}+\sqrt{\mu}}(\x^k-\x^{k-1})$ and $\x^{k+1}=\y^k-\eta\nabla f(\y^k)$.} have the same projection onto the optimum solution set. It does not hold for problem (\ref{non_convex_problem0}).

\subsection{Our Contributions}

In this paper, we use Nesterov's acceleration scheme for problem (\ref{non_convex_problem0}) and an efficient accelerated gradient method with alternating constraint is proposed, which operates on the $\U$ factors directly. We back up our method with provable theoretical results. Specifically, our contributions can be summarized as follows:
\begin{enumerate}
\item We establish the curvature of local restricted strong convexity along a certain trajectory by restricting the problem onto a constraint set, which allows us to use the classical accelerated gradient method for convex programs to solve the constrained problem. We build our result with the tool of polar decomposition.
\item In order to reduce the negative influence of the constraint and ensure the convergence to the critical point of the original unconstrained problem, rather than the reformulated constrained problem, we propose a novel alternating constraint strategy and combine it with the classical accelerated gradient method.
\item When $f$ is restricted $\mu$-strongly convex and restricted $L$-smooth, our method has the local linear convergence to the optimum solution, which has the same dependence on $\sqrt{L/\mu}$ as convex programming. As far as we know, we are the first to establish the convergence matching the optimal dependence on $\sqrt{L/\mu}$ for this kind of nonconvex problems. Globally, our method converges to a critical point of problem (\ref{non_convex_problem0}) from any initializer.
\end{enumerate}

\subsection{Notations and Assumptions}\label{ass_section}
For matrices $\U,\V\in\mathcal{R}^{n\times r}$, we use $\|\U\|_F$ as the Frobenius norm, $\|\U\|_2$ as the spectral norm and $\<\U,\V\>=\tr(\U^T\V)$ as their inner products. We denote $\sigma_r(\U)$ as the smallest singular value of $\U$ and $\sigma_1(\U)=\|\U\|_2$ as the largest one. We use $\U_{ S}\in\mathcal{R}^{r\times r}$ as the submatrix of $\U$ with the rows indicated by the index set $S\subseteq\{1,2,\cdots,n\}$, $\U_{-S}\in\mathcal{R}^{(n-r)\times r}$ as the submatrix with the rows indicated by the indexes out of $S$ and $\X_{ S, S}\in\mathcal{R}^{r\times r}$ as the submatrix of $\X$ with the rows and columns indicated by $ S$. $\X\succeq 0$ means that $\X$ is symmetric and positive semidefinite. For the objective function $g(\U)$, its gradient w.r.t. $\U$ is $\nabla g(\U)=2\nabla f(\U\U^T)\U$. We assume that $\nabla f(\U\U^T)$ is symmetric for simplicity. Our conclusions for the asymmetric case naturally generalize since $\nabla g(\U)=\nabla f(\U\U^T)\U+\nabla f(\U\U^T)^T\U$ in this case. Denote the optimum solution set of problem (\ref{non_convex_problem0}) as
\begin{eqnarray}
\mathcal{X}^*=\{\U^*:\U^*\in\mathcal{R}^{n\times r},\U^*{\U^*}^T=\X^*\}.\label{X_define}
\end{eqnarray}
where $\X^*$ is a minimizer of problem (\ref{convex_problem}). An important issue in minimizing $g(\U)$ is that its optimum solution is not unique, i.e., if $\U^*$ is the optimum solution of problem (\ref{non_convex_problem0}), then $\U^*\R$ is also an optimum solution for any orthogonal matrix $\R\in\mathcal{R}^{r\times r}$. Given $\U$, we define the optimum solution that is closest to $\U$ as
\begin{eqnarray}
P_{\mathcal{X}^*}(\U)=\U^*\R, \mbox{ where }\R=\argmin_{\R\in\mathcal{R}^{r\times r},\R\R^T=\I}\|\U^*\R-\U\|_F^2.\notag
\end{eqnarray}

\subsubsection{Assumptions}
In this paper, we assume that $f$ is restricted $\mu$-strongly convex and $L$-smooth on the set $\{\X:\X\succeq 0,\rank(\X)\leq r\}$. We state the standard definitions below.
\begin{definition}
Let $f:\mathcal{R}^{n\times n}\rightarrow \mathcal{R}$ be a convex differentiable function. Then, $f$ is restricted $\mu$-strongly convex on the set $\{\X:\X\succeq 0,\rank(\X)\leq r\}$ if, for any $\X,\Y\in\{\X:\X\succeq 0,\rank(\X)\leq r\}$, we have
\begin{eqnarray}
&f(\Y)\geq f(\X)+\<\nabla f(\X),\Y-\X\>+\frac{\mu}{2}\|\Y-\X\|_F^2. \notag
\end{eqnarray}
\end{definition}
\begin{definition}
Let $f:\mathcal{R}^{n\times n}\rightarrow \mathcal{R}$ be a convex differentiable function. Then, $f$ is restricted $L$-smooth on the set $\{\X:\X\succeq 0,\rank(\X)\leq r\}$ if, for any $\X,\Y\in\{\X:\X\succeq 0,\rank(\X)\leq r\}$, we have
\begin{eqnarray}
&f(\Y)\leq f(\X)+\<\nabla f(\X),\Y-\X\>+\frac{L}{2}\|\Y-\X\|_F^2 \notag
\end{eqnarray}
and
\begin{eqnarray}
&\|\nabla f(\Y)-\nabla f(\X)\|_F\leq L\|\Y-\X\|_F.\notag
\end{eqnarray}
\end{definition}
\subsubsection{Polar decomposition}
Polar decomposition is a powerful tool for matrix analysis. We briefly review it in this section. We only describe the left polar decomposition of a square matrix.
\begin{definition}
The polar decomposition of a matrix $\A\in\mathcal{R}^{r\times r}$ has the form $\A=\H\Q$ where $\H\in\mathcal{R}^{r\times r}$ is positive semidefinite and $\Q\in\mathcal{R}^{r\times r}$ is an orthogonal matrix.
\end{definition}
If $\A\in\mathcal{R}^{r\times r}$ is of full rank, then $\A$ has the unique polar decomposition with positive definite $\H$. In fact, since a positive semidefinite Hermitian matrix has a unique positive semidefinite square root, $\H$ is uniquely given by $\H=\sqrt{\A\A^T}$. $\Q=\H^{-1}\A$ is also unique.

In this paper, we use the tool of polar decomposition's perturbation theorem to build the restricted strong convexity of $g(\U)$. It is described below.
\begin{lemma}\label{Polar_T}\citep{li-1995}
Let $\A\in\mathcal{R}^{r\times r}$ be of full rank and $\H\Q$ be its unique polar decomposition, $\A+\triangle\A$ be of full rank and $(\H+\triangle\H)(\Q+\triangle\Q)$ be its unique polar decomposition. Then, we have
\begin{eqnarray}
\|\triangle\Q\|_F\leq \frac{2}{\sigma_r(\A)}\|\triangle \A\|_F.\notag
\end{eqnarray}
\end{lemma}

\section{The Restricted Strongly Convex Curvature}\label{curvature_section}
Function $g(\U)$ is a special kind of nonconvex function and the non-convexity only comes from the factorization of $\U\U^T$. Based on this observation, we exploit the special curvature of $g(\U)$ in this section.

The existing works proved the local linear convergence of the \emph{gradient descent method} for problem (\ref{non_convex_problem0}) by exploiting curvatures such as the local second order growth property \cite{sun-2015,Chen-2015} or the $(\alpha,\beta)$ regularity condition \cite{Jin-2017,Srinadh-2016-colt,Bhojanapalli-2016,wang-2016}. The former is described as
\begin{eqnarray}
g(\U)\geq g(\U^*)+\frac{\alpha}{2}\|P_{\mathcal{X}^*}(\U)-\U\|_F^2,\forall \U\label{sccond1}
\end{eqnarray}
while the later is defined as
\begin{eqnarray}
\<\nabla g(\U),\U-P_{\mathcal{X}^*}(\U)\>\geq \frac{\alpha}{2}\|P_{\mathcal{X}^*}(\U)-\U\|_F^2+\frac{1}{2\beta}\|\nabla g(\U)\|_F^2,\forall \U,\label{sccond2}
\end{eqnarray}
where $\U^*\in\mathcal{X}^*$ and $P_{\mathcal{X}^*}(\U)$ is defined in (\ref{X_define}). Both (\ref{sccond1}) and (\ref{sccond2}) can be derived by the local weakly strongly convex condition \cite{Necoara-2016} combing with the smoothness of $g(\U)$. The former is described as
\begin{eqnarray}
g(\U^*)\geq g(\U)+\<\nabla g(\U),P_{\mathcal{X}^*}(\U)-\U\>+\frac{\alpha}{2}\|P_{\mathcal{X}^*}(\U)-\U\|_F^2,\label{strong_C_equ}
\end{eqnarray}
where $\alpha=\mu\sigma^2_r(\U^*)$. As discussed in Section \ref{ass_section}, the optimum solution of problem (\ref{non_convex_problem0}) is not unique. This non-uniqueness makes the difference between the weakly strong convexity and strong convexity, e.g., on the right hand side of (\ref{strong_C_equ}), we use $P_{\mathcal{X}^*}(\U)$, rather than $\U^*$. Moreover, the weakly strongly convex condition cannot infer convexity and $g(\U)$ is not convex even around a small neighborhood of the global optimum solution \citep{Li-2016}.

Necoara, Nesterov and Glineur \cite{Necoara-2016} studied several conditions under which the linear convergence of the \emph{gradient descent method} is guaranteed for general convex programming without strong convexity. The weakly strongly convex condition is the strongest one and can derive all the other conditions. However, it is not enough to analyze the accelerated gradient method only with the weakly strongly convex condition. Necoara et al. \cite{Necoara-2016} proved the acceleration of the classical \emph{accelerated gradient method} under an additional assumption that all the iterates $\{\y^k,k=0,1,\cdots\}$ have the same projection onto the optimum solution set besides the weakly strongly convex condition and the smoothness condition. From the proof in \cite[Section 5.2.1]{Necoara-2016}, we can see that the non-uniqueness of the optimum solution makes the main trouble to analyze the accelerated gradient method\footnote{\cite{Necoara-2016} used induction to prove \cite[Lemma 1]{Necoara-2016}. When the optimum solution is not unique, $\y^*$ in \cite[Equation (57)]{Necoara-2016} should be replaced by $P_{\mathcal{X}^*}(\y^k)$ and they have different values for different $k$. Thus, the induction is not correct any more.}. The additional assumption made in \cite{Necoara-2016} somehow aims to reduce this non-uniqueness. Since this assumption is not satisfied for problem (\ref{non_convex_problem0}), only (\ref{strong_C_equ}) is not enough to prove the acceleration for problem (\ref{non_convex_problem0}) and it requires us to exploit stronger curvature than (\ref{strong_C_equ}) to analyze the accelerated gradient method. 

Motivated by \cite{Necoara-2016}, we should remove the non-uniqueness in problem (\ref{non_convex_problem0}). Our intuition is based on the following observation. Suppose that we can find an index set $ S\subseteq\{1,2,\cdots,n\}$ with size $r$ such that $\X^*_{ S, S}$ is of $r$ full rank, then there exists a unique decomposition $\X^*_{ S, S}=\U^*_{ S}(\U^*_{ S})^T$ where we require $\U^*_{ S}\succ 0$. Thus, we can easily have that there exists a unique $\U^*$ such that $\U^*{\U^*}^T=\X^*$ and $\U^*_{ S}\succ 0$. To verify it, consider $ S=\{1,\cdots,r\}$ for simplicity. Then $\U\U^T=\left(
  \begin{array}{cc}
     \U_{S}\U_{S}^T & \U_{S}\U_{-S}^T\\
     \U_{-S}\U_{S}^T & \U_{-S}\U_{-S}^T
  \end{array}
\right)=\left(
  \begin{array}{cc}
     \X_{S,S} & \X_{S,-S} \\
     \X_{-S,S} & \X_{-S,-S}
  \end{array}
\right)$. The uniqueness of $\U_S$ comes from $\X_{S,S}\succ 0$ and $\U_S\succ 0$ and the uniqueness of $\U_{-S}$ comes from $\U_{-S}=\X_{-S,S}\U_{S}^{-T}$.

Based on the above observation, we can reformulate problem (\ref{non_convex_problem0}) as
\begin{eqnarray}
\min_{\U\in\Omega_{ S}} g(\U)\label{non_convex_problem}
\end{eqnarray}
where
\begin{eqnarray}
\Omega_{ S}=\{\U\in\mathcal{R}^{n\times r}:\U_{ S}\succeq \epsilon \I\}\notag
\end{eqnarray}
and $\epsilon$ is a small enough constant such that $\epsilon\ll \sigma_r(\U_{ S}^*)$. We require $\U_{ S}\succeq \epsilon \I$ rather than $\U_{ S}\succ 0$ to make the projection onto $\Omega_S$ computable. Due to the additional constraint of $\U\in\Omega_S$, we observe that the optimum solution of problem (\ref{non_convex_problem}) is unique. Moreover, the minimizer of (\ref{non_convex_problem}) minimizes also (\ref{non_convex_problem0}).

Until now, we are ready to establish a stronger curvature than (\ref{strong_C_equ}) by restricting the variables of $g(\U)$ on the set $\Omega_S$. We should lower bound $\|P_{\mathcal{X}^*}(\U)-\U\|_F^2$ in (\ref{strong_C_equ}) by $\|\U^*-\U\|_F^2$. Our result is built upon polar decomposition's perturbation theorem \cite{li-1995}. Based on Lemma \ref{Polar_T}, we first establish the following critical lemma.
\begin{lemma}\label{perturbation_T}
For any $\U\in\Omega_{ S}$ and $\V\in\Omega_{ S}$, let $\R=\argmin_{\R\in\mathcal{R}^{r\times r},\R\R^T=\I}\|\V\R-\U\|_F^2$ and $\hat\V=\V\R$. Then, we have
\begin{eqnarray}
\|\V-\U\|_F\leq \frac{3\|\U\|_2}{\sigma_r(\U_S)}\|\hat\V-\U\|_F.\notag
\end{eqnarray}
\end{lemma}
\begin{proof}
Since the conclusion is not affected by permutating the rows of $\U$ and $\V$ under the same permutation, we can consider the case of $ S=\{1,\cdots,r\}$ for simplicity. Let $\U=\left(
  \begin{array}{c}
     \U_1 \\
     \U_2
  \end{array}
\right)$, $\V=\left(
  \begin{array}{c}
     \V_1 \\
     \V_2
  \end{array}
\right)$ and $\hat\V=\left(
  \begin{array}{c}
     \hat\V_1 \\
     \hat\V_2
  \end{array}
\right)$, where $\U_1,\V_1,\hat\V_1\in\mathcal{R}^{r\times r}$. Then, we have $\hat\V_1=\V_1\R$. From $\U\in\Omega_{ S}$ and $\V\in\Omega_{ S}$, we know $\U_1\succ 0$ and $\V_1\succ 0$. Thus, $\U_1\I$ and $\V_1\R$ are the unique polar decompositions of $\U_1$ and $\hat\V_1$, respectively. From Lemma \ref{Polar_T}, we have
\begin{eqnarray}
\|\R-\I\|_F\leq \frac{2}{\sigma_r(\U_1)}\|\hat\V_1-\U_1\|_F.\notag
\end{eqnarray}
With some simple computations, we can have
\begin{eqnarray}
\begin{aligned}\label{cont3}
\|\V-\U\|_F=& \|\hat\V\R^T-\U\|_F\\
=& \|\hat\V\R^T-\U\R^T+\U\R^T-\U\|_F\\
\leq& \|\hat\V\R^T-\U\R^T\|_F+\|\U\R^T-\U\|_F\\
\leq& \|\hat\V-\U\|_F+\|\U\|_2\|\R-\I\|_F\\
\leq& \|\hat\V-\U\|_F+\frac{2\|\U\|_2}{\sigma_r(\U_1)}\|\hat\V_1-\U_1\|_F\\
\leq& \frac{3\|\U\|_2}{\sigma_r(\U_1)}\|\hat\V-\U\|_F,
\end{aligned}
\end{eqnarray}
where we use $\sigma_r(\U_1)\leq \|\U\|_2$ and $\|\hat\V_1-\U_1\|_F\leq\|\hat\V-\U\|_F$ in the last inequality. Replacing $\U_1$ with $\U_{ S}$, we can have the conclusion.
\end{proof}

Built upon Lemma \ref{perturbation_T}, we can give the local restricted strong convexity of $g(\U)$ on the set $\Omega_S$ in the following theorem. There are two differences between the restricted strong convexity and the weakly strong convexity: (i) the restricted strong convexity removes the non-uniqueness and (ii) the restricted strong convexity establishes the curvature between any two points $\U$ and $\V$ in a local neighborhood of $\U^*$, while (\ref{strong_C_equ}) only exploits the curvature between $\U$ and the optimum solution.

\begin{theorem}\label{strong_Cm}
Assume that $\U^*\in\Omega_{ S}\cap\mathcal{X}^*$, $\U\in\Omega_{ S}$ and $\V\in\Omega_{ S}$ with $\|\U-\U^*\|_F\leq C$ and $\|\V-\U^*\|_F\leq C$, where $C=\frac{\mu\sigma_r^2(\U^*)\sigma_r^2(\U_{ S}^*)}{100L\|\U^*\|_2^3}$. Then, we have
\begin{eqnarray}
g(\U)\geq g(\V)+\<\nabla g(\V),\U-\V\>+\frac{\mu\sigma_r^2(\U^*)\sigma_r^2(\U_{ S}^*)}{50\|\U^*\|_2^2}\|\U-\V\|_F^2.\notag
\end{eqnarray}
\end{theorem}
\begin{proof}
From the restricted convexity of $f(\X)$, we have
\begin{eqnarray}
\begin{aligned}\label{cont1}
&f(\V\V^T)-f(\U\U^T)\\
\leq&\<\nabla f(\V\V^T),\V\V^T-\U\U^T\>-\frac{\mu}{2}\|\V\V^T-\U\U^T\|_F^2\\
=&\<\nabla f(\V\V^T),(\V-\U)\V^T\>+\<\nabla f(\V\V^T),\V(\V-\U)^T\>\\
&-\<\nabla f(\V\V^T),(\V-\U)(\V-\U)^T\>-\frac{\mu}{2}\|\V\V^T-\U\U^T\|_F^2\\
=&2\<\nabla f(\V\V^T)\V,\V-\U\>-\<\nabla f(\V\V^T),(\V-\U)(\V-\U)^T\>\\
&-\frac{\mu}{2}\|\V\V^T-\U\U^T\|_F^2\\
\leq&2\<\nabla f(\V\V^T)\V,\V-\U\>-\<\nabla f(\V\V^T)-\nabla f(\X^*),(\V-\U)(\V-\U)^T\>\\
&-\frac{\mu}{2}\|\V\V^T-\U\U^T\|_F^2.
\end{aligned}
\end{eqnarray}
where we use $\nabla f(\X^*)\succeq 0$ proved in Lemma \ref{sdp_lemma} and the fact that the inner product of two positive semidefinite matrices is nonnegative in the last inequality, i.e., $\<\nabla f(\X^*),(\V-\U)(\V-\U)^T\>\geq 0$. Applying Von Neumann’s trace inequality and Lemma \ref{matrix_lemma} to bound the second term, applying Lemmas \ref{perturbation_T} and \ref{diff_lower_bound} to bound the third term, we can have
\begin{eqnarray}
&&f(\V\V^T)-f(\U\U^T)\notag\\
&\leq&2\<\nabla f(\V\V^T)\V,\V-\U\>+L(\|\U^*\|_2+\|\V\|_2)\|\V-\U^*\|_F\|\V-\U\|_F^2\notag\\
&&-\frac{(\sqrt{2}-1)\mu\sigma_r^2(\U)\sigma_r^2(\U_{ S})}{9\|\U\|_2^2}\|\V-\U\|_F^2\notag\\
&\leq&\<\nabla g(\V),\V\hspace*{-0.05cm}-\hspace*{-0.05cm}\U\>\hspace*{-0.05cm}-\hspace*{-0.05cm}\left(\hspace*{-0.05cm}\frac{\mu\sigma_r^2(\U^*)\sigma_r^2(\U_{ S}^*)}{23.1\|\U^*\|_2^2}\hspace*{-0.05cm}-\hspace*{-0.05cm}2.01L\|\U^*\|_2\|\V\hspace*{-0.05cm}-\hspace*{-0.05cm}\U^*\|_F\hspace*{-0.05cm}\right)\|\V\hspace*{-0.05cm}-\hspace*{-0.05cm}\U\|_F^2,\notag
\end{eqnarray}
where we use Lemma \ref{variable_bound} in the last inequality. From the assumption of $\|\V-\U^*\|_F\leq C$, we can have the conclusion. We leave Lemmas \ref{sdp_lemma}, \ref{diff_lower_bound}, \ref{variable_bound} and \ref{matrix_lemma} in Appendix A.
\end{proof}
\subsection{Smoothness of Function $g(\U)$}

Besides the local restricted strong convexity, we can also prove the smoothness of $g(\U)$, which is built in the following theorem.
\begin{theorem}\label{Lipschitz_S}
Let $\hat L=2\|\nabla f(\V\V^T)\|_2+L(\|\V\|_2+\|\U\|_2)^2$. Then, we can have
\begin{eqnarray}
g(\U)\leq g(\V)+\<\nabla g(\V)\V,\U-\V\>+\frac{\hat L}{2}\|\U-\V\|_F^2.\notag
\end{eqnarray}
\end{theorem}
\begin{proof}
From the restricted Lipschitz smoothness of $f$ and a similar induction to (\ref{cont1}), we have
\begin{eqnarray}
&&f(\U\U^T)-f(\V\V^T)\notag\\
&\leq& \<\nabla f(\V\V^T),\U\U^T-\V\V^T\>+\frac{L}{2}\|\U\U^T-\V\V^T\|_F^2\notag\\
&=& \<\nabla f(\V\V^T),(\U-\V)(\U-\V)^T\>\notag\\
&&+2\<\nabla f(\V\V^T)\V,\U-\V\>+\frac{L}{2}\|\U\U^T-\V\V^T\|_F^2.\notag
\end{eqnarray}
Applying Von Neumann's trace inequality to the first term, applying Lemma \ref{matrix_lemma} to the third term, we can have the conclusion.
\end{proof}
When restricted in a small neighborhood of $\U^*$, we can give a better estimate for the smoothness parameter $\hat L$, as follows. The proof is provided in Appendix A.
\begin{corollary}\label{Lipschitz_constant_lemma}
Assume that $\U^*\in\Omega_S\cap\mathcal{X}^*$ and $\U^k,\V^k,\Z^k\in\Omega_{ S}$ with $\|\V^k-\U^*\|_F\leq C$, $\|\U^k-\U^*\|_F\leq C$ and $\|\Z^k-\U^*\|_F\leq C$, where $C$ is defined in Theorem \ref{strong_Cm}. Let $L_g=38L\|\U^*\|_2^2+2\|\nabla f(\X^*)\|_2$ and $\eta=\frac{1}{L_g}$. Then, we have
\begin{subequations}
\begin{align}
&g(\U^{k+1})\leq g(\V^k)+\<\nabla g(\V^k),\U^{k+1}-\V^k\>+\frac{L_g}{2}\|\U^{k+1}-\V^k\|_F^2.\notag
\end{align}
\end{subequations}
\end{corollary}

\section{Accelerated Gradient Method with Alternating Constraint}\label{subsec0}
For problem (\ref{non_convex_problem}), $\Omega_S$ is a convex set and from Theorem \ref{strong_Cm} and Corollary \ref{Lipschitz_constant_lemma} we know that the objective $g(\U)$ behaves locally like a strongly convex and smooth function when restricted on the set $\Omega_S$. Thus, we can use the classical method for convex programming to solve problem (\ref{non_convex_problem}), e.g., the accelerated gradient method\footnote{However, it is still more challenging than convex programming since we should guarantee that all the variables in Theorem \ref{strong_Cm} belong to $\Omega_S$, while it is not required in convex programming. So the conclusion in \cite{Necoara-2016} cannot be applied to problem (\ref{non_convex_problem}) since we cannot obtain $\y^{k+1}\in\Omega_{ S}$ given $\x^{k+1}\in\Omega_{S}$ and $\x^{k}\in\Omega_{S}$ because $\y^{k+1}$ is not a convex combination of $\x^{k+1}$ and $\x^{k}$.}.

However, there remains a practical issue that when solving problem (\ref{non_convex_problem}), we may get stuck at a local minimum of problem (\ref{non_convex_problem}) at the boundary of the constraint $\U\in\Omega_{ S}$, which is not the optimum solution of problem (\ref{non_convex_problem0}). In other words, we may halt before reaching the acceleration region, i.e., the local neighborhood of the optimum solution of problem (\ref{non_convex_problem0}). To overcome this trouble, we propose a novel alternating trajectory strategy. Specifically, we define two sets $\Omega_{S^1}$ and $\Omega_{S^2}$ as follows
\begin{eqnarray}
\Omega_{ S^1}=\{\U\in\mathcal{R}^{n\times r}:\U_{ S^1}\succeq \epsilon \I\},\quad \Omega_{ S^2}=\{\U\in\mathcal{R}^{n\times r}:\U_{ S^2}\succeq \epsilon \I\}\notag
\end{eqnarray}
and minimize the objective $g(\U)$ along the trajectories of $\Omega_{S^1}$ and $\Omega_{S^2}$ alternatively, i.e., when the iteration number $t$ is odd, we minimize $g(\U)$ with the constraint of $\U\in\Omega_{S^1}$, and when $t$ is even, we minimize $g(\U)$ with the constraint of $\U\in\Omega_{S^2}$. Intuitively, when the iterates approach the boundary of $\Omega_{S^1}$, we cancel the constraint of positive definiteness on $\U_{ S^1}$ and put it on $\U_{ S^2}$. Fortunately, with this strategy we can cancel the negative influence of the constraint. We require that both the two index sets $ S^1$ and $ S^2$ are of size $r$ and $ S^1\cap S^2=\emptyset$ such that $\U^*_{ S^1}$ and $\U^*_{ S^2}$ are of full rank. Given proper $S^1$ and $S^2$, we can prove that the method globally converges to a critical point of problem (\ref{non_convex_problem0}). i.e., a point with $\nabla g(\U)=0$, rather than a critical point of problem (\ref{non_convex_problem}).

We describe our method in Algorithm \ref{agd_alg}. We use Nesterov's acceleration scheme in the inner loop with finite $K$ iterations and restart the acceleration scheme at each outer iteration. At the end of each outer iteration, we change the constraint and transform $\U^{t,K+1}\in\Omega_{ S}$ to a new point $\U^{t+1,0}\in\Omega_{ S'}$ via polar decomposition such that $g(\U^{t,K+1})=g(\U^{t+1,0})$. At step (\ref{Z_step}), we need to project $\Z\equiv\Z^{t,k}-\frac{\eta}{\theta_k}\nabla g(\V^{t,k})$ onto $\Omega_S$. Let $\A\Sigma\A^T$ be the eigenvalue decomposition of $\frac{\Z_S+\Z_S^T}{2}$ and $\hat\Sigma=\diag([\max\{\epsilon,\Sigma_{1,1}\},\cdots,\max\{\epsilon,\Sigma_{r,r}\}])$, then $\Z^{t,k+1}_S=\A\hat\Sigma\A^T$ and $\Z^{t,k+1}_{-S}=\Z_{-S}$. At step (\ref{theta}), $\theta_{k+1}$ is computed by $\theta_{k+1}=\frac{\sqrt{\theta_{k}^4+4\theta_{k}^2}-\theta_{k}^2}{2}$. At the end of each outer iteration, we need to compute the polar decomposition. Let $\A\Sigma\B^T$ be the SVD of $\U_{S'}^{t,K+1}$, then we can set $\H=\A\Sigma\A^T$ and $\Q=\A\B^T$. In Algorithm \ref{agd_alg}, we predefine $S^1$ and $S^2$ and fix them during the iterations. In Section \ref{subsec3} we will discuss how to find $S^1$ and $S^2$ using some local information.

At last, let's compare the per-iteration cost of Algorithm \ref{agd_alg} with the methods operating on $\X$ space. Both the eigenvalue decomposition and polar decomposition required in Algorithm \ref{agd_alg} perform on the submatrices of size $r\times r$, which need $O(r^3)$ operations. Thus, the per-iteration complexity of Algorithm \ref{agd_alg} is $O(nr+r^3)$. As a comparison, the methods operating on $\X$ space require at least the top-$r$ singular value/vectors, which need $O(n^2r)$ operations for the deterministic algorithms and $O(n^2\log r)$ for randomized algorithms \cite{Halko2011}. Thus, our method is more efficient at each iteration when $r\ll n$, especially when $r$ is upper bounded by a constant independent on $n$.
\begin{algorithm}[tb]
   \caption{Accelerated Gradient Descent with Alternating Constraint}
   \label{agd_alg}
\begin{algorithmic}
   \STATE Initialize $\Z^{0,0}=\U^{0,0}\in\Omega_{S^2}$, $\eta$, $K$, $\epsilon$.
   \FOR{$t=0,1,2,\cdots$}
   \STATE $\theta_0=1$.
   \FOR{$k=0,1,\dots,K$}
   \STATE
    \begin{eqnarray}
    &&\V^{t,k}=(1-\theta_k)\U^{t,k}+\theta_k\Z^{t,k}.\label{V_step}\\
    &&\Z^{t,k+1}=\argmin_{\Z\in\Omega_{ S}} \hspace*{-0.05cm}\<\hspace*{-0.05cm}\nabla g(\V^{t,k}),\Z\hspace*{-0.05cm}\>\hspace*{-0.05cm}+\hspace*{-0.05cm}\frac{\theta_k}{2\eta}\hspace*{-0.05cm}\left\|\Z\hspace*{-0.05cm}-\hspace*{-0.05cm}\Z^{t,k}\right\|_F^2\hspace*{-0.05cm}, S\hspace*{-0.05cm}=\hspace*{-0.05cm}\left\{\hspace*{-0.05cm}
  \begin{array}{c}
      S^1, \mbox{ if } t\mbox{ is odd}, \\
      S^2, \mbox{ if } t\mbox{ is even}.
  \end{array}
\right.\label{Z_step}\\
    &&\U^{t,k+1}=(1-\theta_k)\U^{t,k}+\theta_k\Z^{t,k+1}.\label{U_step}\\
    &&\mbox{compute }\theta_{k+1}\mbox{ from }\frac{1-\theta_{k+1}}{\theta_{k+1}^2}=\frac{1}{\theta_k^2}.\label{theta}
    \end{eqnarray}
   \ENDFOR
   \STATE Let $\H\Q=\U_{ S'}^{t,K+1}$ be its polar decomposition and $\Z^{t+1,0}=\U^{t+1,0}=\U^{t,K+1}\Q^T$, where $ S'=\left\{
  \begin{array}{c}
      S^2, \mbox{ if } S=S^1, \\
      S^1, \mbox{ if } S=S^2.
  \end{array}
\right.$
   \ENDFOR
\end{algorithmic}
\end{algorithm}

\subsection{Finding the Index Sets $ S^1$ and $ S^2$}\label{subsec3}

In this section, we consider how to find the index sets $ S^1$ and $ S^2$. $S^1\cap S^2=\emptyset$ can be easily satisfied and we only need to ensure that $\U^*_{ S^1}$ and $\U^*_{ S^2}$ are of full rank. Suppose that we have some initializer $\U^0$ close to $\U^*$. We want to use $\U^0$ to find such $ S^1$ and $ S^2$. We first discuss how to select one index set $ S$ based on $\U^0$.  We can use the volume sampling subset selection algorithm \citep{Guruswami-2012,Avron-2013}, which can select $ S$ such that $\sigma_r(\U^0_{ S})\geq \frac{\sigma_r(\U^0)}{\sqrt{2r(n-r+1)}}$ with probability of $1-\delta'$ in $O(nr^3\log(1/\delta'))$ operations. Then, we can bound $\sigma_r(\U^*_{ S})$ in the following lemma since $\U^0$ is close to $\U^*$.

\begin{lemma}\label{index set_theoremm}
If $\|\U^0-\U^*\|_F\leq 0.01\sigma_r(\U^*)$ and $\|\U_S^0-\U_S^*\|_F\leq \frac{0.99\sigma_r(\U^*)}{2\sqrt{2r(n-r+1)}}$, then for the index set $ S$ returned by the volume sampling subset selection algorithm performed on $\U^0$ after $O(nr^3\log(1/\delta'))$ operations, we have $\sigma_r(\U^*_{ S})\geq \frac{0.99\sigma_r(\U^*)}{2\sqrt{2r(n-r+1)}}$ with probability of $1-\delta'$.
\end{lemma}
\begin{proof}
Form Theorem 3.11 in \citep{Avron-2013}, we have $\sigma_r(\U^0_{ S})\geq\frac{\sigma_r(\U^0)}{\sqrt{2r(n-r+1)}}$ with probability of $1-\delta'$ after $O(nr^3\log(1/\delta'))$ operations. So we can obtain
\begin{eqnarray}
\sigma_r(\hspace*{-0.02cm}\U^0_{ S}\hspace*{-0.02cm})\hspace*{-0.06cm}-\hspace*{-0.06cm}\sigma_r(\hspace*{-0.02cm}\U^*_{ S}\hspace*{-0.02cm})\hspace*{-0.06cm}\leq\hspace*{-0.06cm}\|\U_S^0\hspace*{-0.06cm}-\hspace*{-0.06cm}\U_S^*\|_F\hspace*{-0.06cm}\leq\hspace*{-0.06cm} \frac{0.99\sigma_r(\U^*)}{2\hspace*{-0.02cm}\sqrt{2r(n\hspace*{-0.06cm}-\hspace*{-0.06cm}r\hspace*{-0.06cm}+\hspace*{-0.06cm}1)}}\hspace*{-0.06cm}\leq\hspace*{-0.06cm} \frac{\sigma_r(\U^0)}{2\hspace*{-0.02cm}\sqrt{2r(n\hspace*{-0.06cm}-\hspace*{-0.06cm}r\hspace*{-0.06cm}+\hspace*{-0.06cm}1)}}\hspace*{-0.06cm}\leq\hspace*{-0.06cm}\frac{\sigma_r(\U^0_{ S})}{2},\notag
\end{eqnarray}
which leads to
\begin{eqnarray}
\sigma_r(\U^*_{ S})\geq \frac{\sigma_r(\U^0_{ S})}{2}\geq \frac{\sigma_r(\U^0)}{2\sqrt{2r(n-r+1)}}\geq \frac{0.99\sigma_r(\U^*)}{2\sqrt{2r(n-r+1)}},\notag
\end{eqnarray}
where we use $0.99\sigma_r(\U^*)\leq \sigma_r(\U^0)$, which is proved in Lemma \ref{variable_bound} in Appendix A.

\end{proof}

In the column selection problem and its variants, existing algorithms (please see \cite{Avron-2013} and the references therein) can only find one index set. Our purpose is to find both $ S^1$ and $ S^2$. We believe that this is a challenging target in the theoretical computer science community. In our applications, since $n\gg r$, we may expect that the rank of $\U^0_{- S^1}$ is not influenced after dropping $r$ rows from $\U^0$. Thus, we can use the procedure discussed above again to find $ S^2$ from $\U^0_{- S^1}$. From Lemma \ref{index set_theoremm}, we have $\sigma_r(\U^0_{S^1})\geq \frac{\sigma_r(\U^0)}{\sqrt{2r(n-r+1)}}$ and $\sigma_r(\U^0_{S^2})\geq \frac{\sigma_r(\U_{-S^1}^0)}{\sqrt{2r(n-2r+1)}}$. In the asymmetric case, this challenge disappears. Please see the details in Section \ref{asy_sec}. We show in experiments that Algorithm \ref{agd_alg} works well even for the simple choice of $S^1=\{1,\cdots,r\}$ and $S^2=\{r+1,\cdots,2r\}$. The discussion of finding $ S^1$ and $ S^2$ in this section is only for the theoretical purpose.

\subsection{Initialization}

Our theorem ensures the accelerated linear convergence given that the initial point $\U^0\in\Omega_{S^2}$ is within the local neighborhood of the optimum solution, with radius $C$ defined in Theorem \ref{strong_Cm}. We use the initialization strategy in \cite{Srinadh-2016-colt}. Specifically, let $\X^0=\mbox{Project}_{+}\left(\frac{-\nabla f(0)}{\|\nabla f(0)-\nabla f(11^T)\|_F}\right)$ and $\V^0{\V^0}^T$ be the best rank-$r$ approximation of $\X_0$, where $\mbox{Project}_{+}$ means the projection operator onto the semidefinite cone. Then, \cite{Srinadh-2016-colt} proved $\|\V^0-P_{\mathcal{X}^*}(\V^0)\|_F\leq \frac{4\sqrt{2}r\|\U^*\|_2^2}{\sigma_r(\U^*)}\sqrt{\frac{L^2}{\mu^2}-\frac{2\mu}{L}+1}$. Let $\H\Q=\V^0_{S^2}$ be its polar decomposition and $\U^0=\V^0\Q^T$. Then, $\U^0$ belongs to $\Omega_{S^2}$. Although this strategy does not produce an initial point close enough to the target, we show in experiments that our method performs well in practice. It should be noted that for the gradient descent method to solve the general problem (\ref{convex_problem}), the initialization strategy in \cite{Srinadh-2016-colt} also does not satisfy the requirement of the theorems in \cite{Srinadh-2016-colt} for the general objective $f$.

\section{Accelerated Convergence Rate Analysis}\label{convergence_rate_section}

In this section, we prove the local accelerated linear convergence rate of Algorithm \ref{agd_alg}. We first consider the inner loop. It uses the classical accelerated gradient method to solve problem (\ref{non_convex_problem}) with fixed index set $S$ for finite $K$ iterations. Thanks to the stronger curvature built in Theorem \ref{strong_Cm} and the smoothness in Corollary \ref{Lipschitz_constant_lemma}, we can use the standard proof framework to analyze the inner loop, e.g., \cite{Tseng-2008}. Some slight modifications are needed since we should ensure that all the iterates belong to the local neighborhood of $\U^*$. We present the result in the following lemma and give its proof sketch. For simplicity, we omit the outer iteration number $t$.

\begin{lemma}\label{conver_rate_square}
Assume that $\U^*\in\Omega_{ S}\cap\mathcal{X}^*$ and $\U^0\in\Omega_{ S}$ with $\epsilon\leq0.99\sigma_r(\U_{ S'}^*)$ and $\|\U^0-\U^*\|_F\leq C$. Let $\eta=\frac{1}{L_g}$, where $C$ is defined in Theorem \ref{strong_Cm} and $L_g$ is defined in Corollary \ref{Lipschitz_constant_lemma}. Then, we have $\sigma_r(\U_{ S'}^{K+1})\geq \epsilon$, $\|\U^{K+1}-\U^*\|_F\leq C$ and
\begin{eqnarray}
g(\U^{K+1})-g(\U^*)\leq \frac{2}{(K+1)^2\eta}\left\|\U^*-\U^0\right\|_F^2.\notag
\end{eqnarray}
\end{lemma}
\begin{proof}
We follow four step to prove the lemma.

Step 1: We can easily check that if $\U^0\in\Omega_{ S}$, then all the iterates of $\{\U^k,\V^k,\Z^k\}$ belong to $\Omega_{ S}$ by $0\leq \theta_k\leq 1$, the convexity of $\Omega_S$ and the convex combinations in (\ref{V_step}) and (\ref{U_step}).

Step 2: Consider the $k$-th iteration. If $\|\V^k-\U^*\|_F\leq C$, $\|\Z^k-\U^*\|_F\leq C$ and $\|\U^k-\U^*\|_F\leq C$, then Theorem \ref{strong_Cm} and Corollary \ref{Lipschitz_constant_lemma} hold. From the standard analysis of the accelerated gradient method for convex programming, e.g., Proposition 1 in \cite{Tseng-2008}, we can have
\begin{eqnarray}
\begin{aligned}\label{cont2}
&\frac{1}{\theta_k^2}\left(g(\U^{k+1})-g(\U^*)\right)+\frac{1}{2\eta}\|\Z^{k+1}-\U^*\|_F^2\\
\leq&\frac{1}{\theta_{k-1}^2}\left(g(\U^k)-g(\U^*)\right)+\frac{1}{2\eta}\|\U^*-\Z^k\|_F^2.
\end{aligned}
\end{eqnarray}

Step 3: Since Theorem \ref{strong_Cm} and Corollary \ref{Lipschitz_constant_lemma} hold only in a local neighbourhood of $\U^*$, we need to check that $\{\U^k,\V^k,\Z^k\}$ belongs to this neighborhood for all the iterations, which can be easily done via induction. In fact, from (\ref{cont2}) and the convexity combinations in (\ref{V_step}) and (\ref{U_step}), we know that if the following conditions hold,
\begin{eqnarray}
&&\|\V^k-\U^*\|_F\leq C,\quad \|\U^k-\U^*\|_F\leq C,\quad \|\Z^k-\U^*\|_F\leq C,\notag\\
&&\frac{1}{\theta_{k-1}^2}\left(g(\U^k)-g(\U^*)\right)+\frac{1}{2\eta}\|\Z^k-\U^*\|_F^2\leq \frac{C^2}{2\eta},\notag
\end{eqnarray}
then we can have
\begin{eqnarray}
&&\|\V^{k+1}-\U^*\|_F\leq C,\quad \|\U^{k+1}-\U^*\|_F\leq C,\quad \|\Z^{k+1}-\U^*\|_F\leq C,\notag\\
&&\frac{1}{\theta_k^2}\left(g(\U^{k+1})-g(\U^*)\right)+\frac{1}{2\eta}\|\Z^{k+1}-\U^*\|_F^2\leq \frac{C^2}{2\eta}.\notag
\end{eqnarray}

Step 4: From $\frac{1}{\theta_{-1}}=0$ and Step 3, we know (\ref{cont2}) holds for all the iterations. Thus, we have
\begin{eqnarray}
g(\U^{K+1})-g(\U^*)\leq \frac{\theta_K^2}{2\eta}\left\|\Z^0-\U^*\right\|_F^2\leq \frac{2}{(K+1)^2\eta}\|\Z^0-\U^*\|_F^2,\notag
\end{eqnarray}
where we use $\theta_k\leq\frac{2}{k+1}$ from $\frac{1-\theta_{k+1}}{\theta_{k+1}^2}=\frac{1}{\theta_k^2}$ and $\theta_0=1$.

On the other hand, from the perturbation theorem of singular values, we have
\begin{eqnarray}
\sigma_r(\hspace*{-0.02cm}\U_{ S'}^*\hspace*{-0.02cm})\hspace*{-0.05cm}-\hspace*{-0.05cm}\sigma_r(\hspace*{-0.02cm}\U_{ S'}^{K+1}\hspace*{-0.02cm})\leq\|\U_{ S'}^{K+1}\hspace*{-0.05cm}-\hspace*{-0.05cm}\U_{ S'}^*\|_F\leq\|\U^{K+1}\hspace*{-0.05cm}-\hspace*{-0.05cm}\U^*\|_F\leq C\leq 0.01\sigma_r(\hspace*{-0.02cm}\U_{ S'}^*\hspace*{-0.02cm}),\notag
\end{eqnarray}
which leads to $\sigma_r(\U_{ S'}^{K+1})\geq 0.99\sigma_r(\U_{ S'}^*)\geq\epsilon$.

\end{proof}

Now we consider the outer loop of Algorithm \ref{agd_alg}. Based on Lemma \ref{conver_rate_square}, the second order growth property (\ref{sccond1}) and the perturbation theory of polar decomposition, we can establish the exponentially decreasing of $\|\U^{t,0}-\U^{t,*}\|_F$ in the following lemma.
\begin{lemma}\label{conver_rate_square1}
Assume that $\U^{t,*}\in\Omega_{ S}\cap\mathcal{X}^*$, $\U^{t+1,*}\in\Omega_{ S'}\cap\mathcal{X}^*$ and $\U^{t,0}\in\Omega_{ S}$ with $\epsilon\leq0.99\sigma_r(\U_{ S'}^{t,*})$ and $\|\U^{t,0}-\U^{t,*}\|_F\leq C$. Let $K+1=\frac{28\|\U^{*}\|_2}{\sqrt{\eta\mu}\sigma_r(\U^{*})\min\{\sigma_r(\U_{ S^1}^{*}),\sigma_r(\U_{ S^2}^{*})\}}$. Then, we can have $\U^{t+1,0}\in\Omega_{ S'}$ and
\begin{eqnarray}
\|\U^{t+1,0}-\U^{t+1,*}\|_F\leq \frac{1}{4}\|\U^{t,0}-\U^{t,*}\|_F.\label{aaa3}
\end{eqnarray}
\end{lemma}
\begin{proof}
We follow four steps to prove the lemma.

Step 1. From Lemma \ref{conver_rate_square}, we have $\sigma_r(\U_{ S'}^{t,K+1})\geq \epsilon$, $\|\U^{t,K+1}-\U^{t,*}\|_F\leq C$ and
\begin{eqnarray}
g(\U^{t,K+1})-g(\U^{t,*})\leq \frac{2}{(K+1)^2\eta}\|\U^{t,0}-\U^{t,*}\|_F^2.\label{aaa1}
\end{eqnarray}
From Algorithm \ref{agd_alg}, we have $\sigma_r(\U_{ S'}^{t+1,0})=\sigma_r(\U_{ S'}^{t,K+1})$. So $\U_{ S'}^{t+1,0}\succeq \epsilon \I$ and $\U^{t+1,0}\in\Omega_{ S'}$.

Step 2. From Lemma \ref{strong_CW} in Appendix B, we have
\begin{eqnarray}
g(\U^{t,K+1})-g(\U^{t,*})\geq 0.4\mu\sigma_r^2(\U^{t,*})\|\U^{t,K+1}-\hat\U^{t,*}\|_F^2,\label{aaa2}
\end{eqnarray}
where $\hat\U^{t,*}=P_{\mathcal{X}^*}(\U^{t,K+1})=\U^{t,*}\R$ and $\R=\argmin_{\R\R^T=\I}\|\U^{t,*}\R-\U^{t,K+1}\|_F^2$.

Step 3. Given (\ref{aaa1}) and (\ref{aaa2}), in order to prove (\ref{aaa3}), we only need to lower bound $\|\U^{t,K+1}-\hat\U^{t,*}\|_F$ by $\|\U^{t+1,0}-\U^{t+1,*}\|_F$.

From Algorithm \ref{agd_alg}, we know that $\H\Q=\U_{ S'}^{t,K+1}$ is the unique polar decomposition of $\U_{ S'}^{t,K+1}$ and $\U^{t+1,0}=\U^{t,K+1}\Q^T$. Let $\H^*\Q^*=\hat\U_{ S'}^{t,*}$ be its unique polar decomposition and $\U^{t+1,*}=\hat\U^{t,*}(\Q^*)^T$, then $\U^{t+1,*}\in\Omega_{ S'}\cap\mathcal{X}^*$. From the perturbation theorem of polar decomposition in Lemma \ref{Polar_T}, we have
\begin{eqnarray}
\|\Q-\Q^*\|_F\leq \frac{2}{\sigma_r(\hat\U_{ S'}^{t,*})}\|\U_{ S'}^{t,K+1}-\hat\U_{ S'}^{t,*}\|_F.\notag
\end{eqnarray}
Similar to (\ref{cont3}), we have
\begin{eqnarray}
\begin{aligned}\label{cont4}
&\|\U^{t+1,0}-\U^{t+1,*}\|_F\\
=&\|\U^{t,K+1}\Q^T-\hat\U^{t,*}(\Q^*)^T\|_F\\
=&\|\U^{t,K+1}\Q^T-\hat\U^{t,*}\Q^T+\hat\U^{t,*}\Q^T-\hat\U^{t,*}(\Q^*)^T\|_F\\
\leq&\|\U^{t,K+1}-\hat\U^{t,*}\|_F+\|\hat\U^{t,*}\|_2\|\Q-\Q^*\|_F\\
\leq&\frac{3\|\U^{t,*}\|_2}{\sigma_r(\U^{t,*}_{ S'})}\|\U^{t,K+1}-\hat\U^{t,*}\|_F.
\end{aligned}
\end{eqnarray}

Step 4. Combining (\ref{aaa1}), (\ref{aaa2}) and (\ref{cont4}), we have
\begin{eqnarray}
&&\|\U^{t+1,0}-\U^{t+1,*}\|_F\notag\\
&\leq& \frac{3\|\U^{t,*}\|_2}{\sigma_r(\U^{t,*}_{ S'})}\|\U^{t,K+1}-\hat\U^{t,*}\|_F\notag\\
&\leq& \frac{3\|\U^{t,*}\|_2}{\sigma_r(\U^{t,*}_{ S'})}\frac{\sqrt{5}}{\sqrt{\eta\mu}(K+1)\sigma_r(\U^{t,*})}\|\U^{t,0}-\U^{t,*}\|_F\notag\\
&\leq&\frac{7\|\U^{t,*}\|_2}{\sqrt{\eta\mu}(K+1)\sigma_r(\U^{t,*})\min\{\sigma_r(\U_{ S}^{t,*}),\sigma_r(\U_{ S'}^{t,*})\}}\|\U^{t,0}-\U^{t,*}\|_F.\notag
\end{eqnarray}
From the setting of $K+1$, we can have the conclusion.
\end{proof}

Combing Lemmas \ref{conver_rate_square} and \ref{conver_rate_square1}, we can give the accelerated convergence rate in the following theorem, i.e., after $O\left(\sqrt{\frac{L}{\mu}}\frac{\|\U^*\|_2^2}{\sigma_r(\U^*)\min\{\sigma_r(\U_{ S^1}^*),\sigma_r(\U_{ S^2}^*)\}}\log\frac{1}{\varepsilon}\right)$ total inner iterations, Algorithm \ref{agd_alg} finds an $\varepsilon$-optimal solution in the sense of $g(\U^{t+1,0})-g(\U^{*})\leq\varepsilon$. The proof is provided in Appendix B. In Algorithm \ref{agd_alg}, when the outer iteration number $t$ is odd, the iterates $\{\U^{t,k}\}$ converge to the unique optimum solution of $\Omega_{S^1}\cap\mathcal{X}^*$. When $t$ is even, $\{\U^{t,k}\}$ converge to another optimum solution of $\Omega_{S^2}\cap\mathcal{X}^*$. In our algorithm, we set $\eta$ and $K$ based on a reliable knowledge on $\|\U^*\|_2$ and $\sigma_r(\U^*)$. As suggested by \cite{Srinadh-2016-colt,Park-2016-siam}, they can be estimated by $\|\U^0\|_2$ and $\sigma_r(\U^0)$-up to constants-since $\U^0$ is close to $\U^*$.

\begin{theorem}\label{conver_ratem}
Assume that $\U^*\in\Omega_{S^2}\cap\mathcal{X}^*$ and $\U^{0,0}\in\Omega_{S^2}$ with $\|\U^{0,0}-\U^*\|_F\leq C$ and $\epsilon\leq\min\{0.99\sigma_r(\U^*_{ S^1}),0.99\sigma_r(\U^*_{ S^2})\}$. Then, we have
\begin{eqnarray}
\|\U^{t+1,0}-\U^{t+1,*}\|_F\leq \left(1-\frac{1}{6}\sqrt{\frac{\mu_g}{L_g}}\right)^{(t+1)(K+1)}\|\U^{0,0}-\U^*\|_F,\notag
\end{eqnarray}
and
\begin{eqnarray}
g(\U^{t+1,0})-g(\U^{*})\leq L_g\left(1-\frac{1}{6}\sqrt{\frac{\mu_g}{L_g}}\right)^{2(t+1)(K+1)}\|\U^{0,0}-\U^*\|_F^2,\notag
\end{eqnarray}
where $\mu_g=\frac{\mu\sigma_r^2(\U^*)\min\{\sigma_r^2(\U_{ S^1}^*),\sigma_r^2(\U_{ S^2}^*)\}}{25\|\U^*\|_2^2}$, $L_g=38L\|\U^*\|_2^2+2\|\nabla f(\X^*)\|_2$, $\U^{t,*}=\Omega_{S^1}\cap\mathcal{X}^*$ when $t$ is odd and $\U^{t,*}=\Omega_{S^2}\cap\mathcal{X}^*$ when $t$ is even.
\end{theorem}

\subsection{Comparison to the Gradient Descent}
Bhojanapalli et al. \cite{Srinadh-2016-colt} used the gradient descent to solve problem (\ref{non_convex_problem0}), which consists of the following recursion:
\begin{eqnarray}
\U^{k+1}=\U^k-\eta\nabla g(\U^k).\notag
\end{eqnarray}
With the restricted strong convexity and smoothness of $f(\X)$, Bhojanapalli et al. \cite{Srinadh-2016-colt} proved the linear convergence of gradient descent in the form of
\begin{eqnarray}
\begin{aligned}
&\|\U^{N+1}-P_{\mathcal{X}^*}(\U^{N+1})\|_F^2\\
\leq& \left(1-\frac{\sigma_r^2(\U^*)}{\|\U^*\|_2^2}\frac{\mu}{L+\|\nabla f(\X^*)\|_2/\|\U^*\|_2^2}\right)^N\|\U^0-P_{\mathcal{X}^*}(\U^0)\|_F^2.\label{gd_rate}
\end{aligned}
\end{eqnarray}

As a comparison, from Theorem \ref{conver_ratem}, our method converges linearly within the error of

\noindent$\left(1-\frac{\sigma_r(\U^*)\min\{\sigma_r(\U_{ S^1}^*),\sigma_r(\U_{ S^2}^*)\}}{\|\U^*\|_2^2}\sqrt{\frac{\mu}{L+\|\nabla f(\X^*)\|_2/\|\U^*\|_2^2}}\right)^N$, where $N$ is the total number of inner iterations. From Lemma \ref{index set_theoremm}, we know $\sigma_r(\U^*_{ S})\approx \frac{1}{\sqrt{rn}}\sigma_r(\U^*)$ in the worst case and it is tight \cite{Avron-2013}. Thus, our method has the convergence rate of $\left(1-\frac{\sigma_r^2(\U^*)}{\|\U^*\|_2^2}\sqrt{\frac{\mu}{nr(L++\|\nabla f(\X^*)\|_2/\|\U^*\|_2^2)}}\right)^N$ in the worst case. When the function $f$ is ill-conditioned, i.e., $\frac{L}{\mu}\geq nr$, our method outperforms the gradient descent. This phenomenon is similar to the case observed in the stochastic optimization community: the non-accelerated methods such as SDCA \cite{zhang-2013-jmlr}, SVRG \cite{zhang-2014-siam} and SAG \cite{Schmidt-2017-mp} have the complexity of $O\left(\frac{L}{\mu}\log\frac{1}{\epsilon}\right)$ while the accelerated methods such as Accelerated SDCA \cite{zhang-2015-MP}, Catalyst \cite{lin-2015-nips} and Katyusha \cite{zhu-2017-stoc} have the complexity of $O\left(\sqrt{\frac{mL}{\mu}}\log\frac{1}{\epsilon}\right)$, where $m$ is the sample size. The latter is tight when $\frac{L}{\mu}\geq m$ for stochastic programming \cite{Woodworth-2016}. In matrix completion, the optimal sample complexity is $O(rn\log n)$ \cite{Candes-2009-matrix}. It is unclear whether our convergence rate for problem (\ref{non_convex_problem0}) is tight or there exists a faster method. We leave it as an open problem.

For better reference, we summarize the comparisons in Table \ref{table-com}. We can see that our method has the same optimal dependence on $\sqrt{\frac{L}{\mu}}$ as convex programming.

\begin{table}
\caption{Convergence rate comparisons of the gradient descent method (GD) and accelerated gradient descent method (AGD).\label{table-com}}
\begin{center}
\begin{tabular}{c||c|c}\hline
Method & Convex problem& Nonconvex problem (\ref{non_convex_problem0})
\\\hline\hline
 GD   &  $\left(\frac{L-\mu}{L+\mu}\right)^N$\cite{Nesterov-2004}    &  $\left(1-\frac{\sigma_r^2(\U^*)}{\|\U^*\|_2^2}\frac{\mu}{L+\|\nabla f(\X^*)\|_2/\|\U^*\|_2^2}\right)^N$\cite{Srinadh-2016-colt} \\
\hline
 \multirow{3}{*}{AGD}  &  \multirow{3}{*}{$\left(1-\sqrt{\frac{\mu}{L}}\right)^N$\cite{Nesterov-2004}} &  $\left(1-\frac{\sigma_r(\U^*)\min\{\sigma_r(\U_{ S^1}^*),\sigma_r(\U_{ S^2}^*)\}}{\|\U^*\|_2^2}\sqrt{\frac{\mu}{L+\|\nabla f(\X^*)\|_2/\|\U^*\|_2^2}}\right)^N$\\
                       &&   \hspace*{-2cm}$=\left(1-\frac{\sigma_r^2(\U^*)}{\|\U^*\|_2^2}\sqrt{\frac{\mu}{nr(L+\|\nabla f(\X^*)\|_2/\|\U^*\|_2^2)}}\right)^N$\\
\hline
\end{tabular}
\end{center}
\end{table}
\subsubsection{Dropping the Dependence on $n$}
Our convergence rate has an additional dependence on $n$ compared with the gradient descent method. It comes from $\sigma_r(\U^*_S)$, i.e., Lemma \ref{perturbation_T}. In fact, we use a loose relaxation in the last inequality of (\ref{cont3}), i.e., $\frac{2\|\U\|_2}{\sigma_r(\U_S)}\|\hat\V_S-\U_S\|_F\leq \frac{2\|\U\|_2}{\sigma_r(\U_S)}\|\hat\V-\U\|_F$. Since $\U_S\in\mathcal{R}^{r\times r}$ and $\U\in\mathcal{R}^{n\times r}$, a more suitable estimation should be
\begin{eqnarray}
\frac{2\|\U\|_2}{\sigma_r(\U_S)}\|\hat\V_S-\U_S\|_F\approx \frac{2\|\U\|_2}{\sigma_r(\U_S)}\sqrt{\frac{r}{n}}\|\hat\V-\U\|_F\approx \frac{2r\|\U\|_2}{\sigma_r(\U)}\|\hat\V-\U\|_F.\label{cont15}
\end{eqnarray}
In practice, (\ref{cont15}) holds when the entries of $\U^{t,k}$ and $\V^{t,k}$ converge nearly equally fast to those of $\U^{t,*}$, which may be expected in practice. Thus, under the condition of (\ref{cont15}), our convergence rate can be improved to
\begin{eqnarray}
\left(1-\frac{\sigma_r^2(\U^*)}{r\|\U^*\|_2^2}\sqrt{\frac{\mu}{L+\|\nabla f(\X^*)\|_2/\|\U^*\|_2^2}}\right)^N. \notag
\end{eqnarray}
We numerically verify (\ref{cont15}) in Section \ref{section_exp_4}.

\subsubsection{Examples with Ill-conditioned Objective $f$}

Although the condition number $\frac{L}{\mu}$ approximate to 1 for some famous problems in machine learning, e.g., matrix regression and matrix completion \cite{Chen-2015}, we can still find many problems with ill-conditioned objective, especially in the computer vision applications. We give the example of low rank representation (LRR) \cite{lin-lrr}. The LRR model is a famous model in computer vision. It can be formulated as
\begin{eqnarray}
\min_{\X} \rank(\X) \quad s.t.\quad \D\X=\A,\notag
\end{eqnarray}
where $\A$ is the observed data and $\D$ is a dictionary that linearly spans the data space. We can reformulate the problem as follows:
\begin{eqnarray}
\min_{\X} \|\D\X-\A\|_F^2 \quad s.t.\quad \rank(\X)\leq r.\notag
\end{eqnarray}
We know $L/\mu=\kappa(\D^T\D)$, which is the condition number of $\D^T\D$. If we generate $\D\in\mathcal{R}^{n\times n}$ as a random matrix with normal distribution, then $E\left[\mbox{log}\kappa(\D)\right]\sim\mbox{log}n$ as $n\rightarrow\infty$ \cite{Edelman-1988} and thus $E\left[\frac{L}{\mu}\right]\sim n^2$. We numerically verify on MATLAB that if $n=1000$, then $\frac{L}{\mu}$ is of the order $10^7$, which is much larger than $O(n)$.

\section{Global Convergence}\label{global_section}

In this section, we study the global convergence of Algorithm \ref{agd_alg} without the assumption that $f(\X)$ is restricted strongly convex. We allow the algorithm to start from any initializer. Since we have no information about $\U^*$ when $\U^0$ is far from $\U^*$, we use an adaptive index sets selection procedure for Algorithm \ref{agd_alg}. That is to say, after each inner loop, we check whether $\sigma_r(\U^{t,K+1}_{S'})<\epsilon$ holds. If not, we select the new index set $S'$ using the volume sampling subset selection algorithm.

We first consider the inner loop and establish Lemma \ref{global_lemma}. We drop the outer iteration number $t$ for simplicity and leave the proof in Appendix C.

\begin{lemma}\label{global_lemma}
Assume that $\{\U^k,\V^k\}$ is bounded and $\U^0\in\Omega_{ S}$. Let $\eta\leq\frac{1-\beta_{\max}^2}{\hat L(2\beta_{\max}+1)+2\gamma}$, where $\hat L=2D+4LM^2$, $D=\max\{\|\nabla f(\U^k(\U^k)^T)\|_2,\|\nabla f(\V^k(\V^k)^T)\|_2,\forall k\}$, $M=\max\{\|\U^k\|_2,\|\V^k\|_2,\forall k\}$, $\beta_{\max}=\max\left\{\beta_k,k=0,\cdots,K\right\}$, $\beta_k=\frac{\theta_k(1-\theta_{k-1})}{\theta_{k-1}}$ and $\gamma$ is a small constant. Then, we have
\begin{eqnarray}
g(\U^{K+1})-g(\U^0)\leq-\sum_{k=0}^K \gamma\|\U^{k+1}-\U^k\|_F^2.\notag
\end{eqnarray}
\end{lemma}

Now we consider the outer loop. As discussed in Section \ref{subsec0}, when solving problem (\ref{non_convex_problem}) directly, we may get stuck at the boundary of the constraint. Thanks to the alternating constraint strategy, we can cancel the negative influence of the constraint and establish the global convergence to a critical point of problem (\ref{non_convex_problem0}), which is described in Theorem \ref{global_convergencem}. It establishes that after at most $O\left(\frac{1}{\varepsilon^2}\log\frac{1}{\varepsilon}\right)$ operations, $\U^{T,K+1}$ is an approximate zero gradient point in the precision of $\varepsilon$.

\begin{theorem}\label{global_convergencem}
Assume that $\{\U^{t,k},\V^{t,k}\}$ is bounded and $\sigma_r(\U_{S'}^{t,K+1})\geq\epsilon,\forall t$. Let $\eta$ be the one defined in Lemma \ref{global_lemma}. Then, after at most $T=2\frac{f(\U^{t,0}(\U^{t,0})^T)-f(\X^*)}{\varepsilon^2}$ outer iterations, we have
\begin{eqnarray}
\left\|\nabla g(\U^{T,K+1})\right\|_F\leq 21\left(\frac{1}{\eta\theta_K}+\frac{\hat L}{\theta_K}\right)\varepsilon\notag
\end{eqnarray}
with probability of $1-\delta$. The volume sampling subset selection algorithm needs $O\left(nr^3\log\left(\frac{f(\U^{t,0}(\U^{t,0})^T)-f(\X^*))}{\delta\varepsilon^2}\right)\right)$ operations for each running.
\end{theorem}
\begin{proof}
We follow three steps to prove the theorem.

Step 1. Firstly, we bound the difference of two consecutive variables, i.e., $\U^{t,k+1}-\U^{t,k}$.

From Lemma \ref{global_lemma} we have
\begin{eqnarray}
\gamma\sum_{k=0}^K \|\U^{t,k+1}-\U^{t,k}\|_F^2\leq g(\U^{t,0})-g(\U^{t,K+1}).\notag
\end{eqnarray}
Summing over $t=0,\cdots,T$ yields
\begin{eqnarray}
&&\gamma\sum_{t=0}^{T}\sum_{k=0}^K \|\U^{t,k+1}-\U^{t,k}\|_F^2\leq\sum_{t=0}^{T}\left(g(\U^{t,0})-g(\U^{t,K+1})\right)\notag\\
&&=\sum_{t=0}^{T}\left(g(\U^{t,0})-g(\U^{t+1,0})\right)\leq g(\U^{t,0})-f(\U^*{\U^*}^T).\notag
\end{eqnarray}
So after $T=2\frac{g(\U^{t,0})-f(\X^*)}{\varepsilon^2}$ outer iterations, we must have
\begin{eqnarray}
\sum_{k=0}^K \|\U^{t,k+1}-\U^{t,k}\|_F^2+\sum_{k=0}^K \|\U^{t+1,k+1}-\U^{t+1,k}\|_F^2\leq\varepsilon^2\label{cont7}
\end{eqnarray}
for some $t<T$. Thus, we can bound $\|\U^{t',k+1}-\U^{t',k}\|_F$ by $\varepsilon$, where $t'=t$ or $t'=t+1$. Moreover, from Lemma \ref{temp_lemma2} in Appendix C, we can bound $\|\U^{t',k+1}-\Z^{t',k+1}\|_F$, $\|\Z^{t',k+1}-\Z^{t',k}\|_F$ and $\|\Z^{t',k+1}-\V^{t',k}\|_F$ by $\frac{\varepsilon}{\theta_k}$.

Step 2. Secondly, we bound parts of elements of the gradient, i.e., $\left(\nabla g(\Z^{t,K+1})\right)_{- S^1}$ and $\left(\nabla g(\Z^{t,K+1})\right)_{- S^2}$.

From the optimality condition of (\ref{Z_step}), we have
\begin{eqnarray}
-\frac{\theta_k}{\eta}\hspace*{-0.05cm}\left(\hspace*{-0.05cm}\Z^{t',k+1}\hspace*{-0.05cm}-\hspace*{-0.05cm}\Z^{t',k}\hspace*{-0.05cm}\right)\hspace*{-0.05cm}+\hspace*{-0.05cm}\nabla g(\Z^{t',k+1})\hspace*{-0.05cm}-\hspace*{-0.05cm}\nabla g(\V^{t',k})\hspace*{-0.05cm}\in\hspace*{-0.05cm} \nabla g(\Z^{t',k+1})\hspace*{-0.05cm}+\hspace*{-0.05cm}\partial I_{\Omega_{ S^j}}(\Z^{t',k+1})\notag
\end{eqnarray}
for $j=1$ when $t'=t$ and $j=2$ when $t'=t+1$. From Lemmas \ref{matrix_lemma} and \ref{temp_lemma2}, we can easily check that
\begin{eqnarray}
\left\|-\frac{\theta_k}{\eta}\left(\Z^{t',k+1}-\Z^{t',k}\right)+\nabla g(\Z^{t',k+1})-\nabla g(\V^{t',k})\right\|_F\leq \frac{14\varepsilon}{\eta\theta_k}.\notag
\end{eqnarray}
Thus, we obtain
\begin{eqnarray}
\mbox{dist}\left(0,\nabla g(\Z^{t',k+1})+\partial I_{\Omega_{ S^i}}(\Z^{t',k+1})\right)\leq \frac{14\varepsilon}{\eta\theta_k},\forall k=0,\cdots,K,\notag
\end{eqnarray}
which leads to
\begin{eqnarray}
\left\|\left(\nabla g(\Z^{t,K+1})\right)_{- S^1}\right\|_F\leq \frac{14\varepsilon}{\eta\theta_K},\label{cp1}
\end{eqnarray}
and
\begin{eqnarray}
\left\|\left(\nabla g(\Z^{t+1,1})\right)_{- S^2}\right\|_F\leq \frac{14\varepsilon}{\eta\theta_K},\label{cont8}
\end{eqnarray}
where $\A_{- S}$ means the submatrix with the rows indicated by the indexes out of $ S$. On the other hand,
\begin{eqnarray}
\begin{aligned}\label{cont9}
&\left\|\left(\nabla g(\Z^{t+1,0})\right)_{- S^2}\right\|_F-\left\|\left(\nabla g(\Z^{t+1,1})\right)_{- S^2}\right\|_F\\
\leq&\left\|\left(\nabla g(\Z^{t+1,0})-\nabla g(\Z^{t+1,1})\right)_{- S^2}\right\|_F\\
\leq&\left\|\nabla g(\Z^{t+1,0})-\nabla g(\Z^{t+1,1})\right\|_F\leq\hat L\|\Z^{t+1,0}-\Z^{t+1,1}\|_F\leq\frac{5\hat L\varepsilon}{\theta_K},
\end{aligned}
\end{eqnarray}
where we use Lemma \ref{temp_lemma2} in the last inequality. Combing (\ref{cont8}) and (\ref{cont9}), we can obtain
\begin{eqnarray}
\left\|\left(\nabla g(\Z^{t+1,0})\right)_{- S^2}\right\|_F\leq \frac{19\varepsilon}{\eta\theta_K}.\notag
\end{eqnarray}
Since $\Z^{t+1,0}=\Z^{t,K+1}\Q^T$ for some orthogonal $\Q$, we can have
\begin{eqnarray}
\begin{aligned}\label{cp2}
\frac{19\varepsilon}{\eta\theta_K}\geq&\left\|\left(\nabla g(\Z^{t+1,0})\right)_{- S^2}\right\|_F=\left\|\left(\nabla g(\Z^{t,K+1})\Q^T\right)_{- S^2}\right\|_F\\
=&\left\|\left(\nabla g(\Z^{t,K+1})\right)_{- S^2}\Q^T\right\|_F=\left\|\left(\nabla g(\Z^{t,K+1})\right)_{- S^2}\right\|_F.
\end{aligned}
\end{eqnarray}

Step 3. We bound all the elements of the gradient. Recall that we require $S^1\cap S^2=\emptyset$. Thus, we have $- S^1\cup- S^2=\{1,2,\cdots,n\}$. Then, from (\ref{cp1}) and (\ref{cp2}), we have
\begin{eqnarray}
\left\|\nabla g(\Z^{t,K+1})\right\|_F\leq \left\|\left(\nabla g(\Z^{t,K+1})\right)_{- S^1}\right\|_F+\left\|\left(\nabla g(\Z^{t,K+1})\right)_{- S^2}\right\|_F\leq \frac{33\varepsilon}{\eta\theta_K}.\notag
\end{eqnarray}
At last, we can bound $\left\|\nabla g(\U^{t,K+1})\right\|_F$ from Lemmas \ref{matrix_lemma} and \ref{temp_lemma2}.

From the Algorithm, we know that the index set is selected at most $T$ times. The volume sampling subset selection algorithm succeeds with the probability of $1-\delta'$. So the Algorithm succeeds with the probability at least of $1-T\delta'=1-\delta$. On the other hand, the volume sampling subset selection algorithm needs $O\left(nr^3\log\left(\frac{1}{\delta'}\right)\right)=O\left(nr^3\log\left(\frac{T}{\delta}\right)\right)=O\left(nr^3\log\left(\frac{f(\U^{t,0}(\U^{t,0})^T)-f(\U^*{\U^*}^T))}{\delta\varepsilon^2}\right)\right)$ operations.
\end{proof}

\section{Minimizing (\ref{non_convex_problem0}) Directly without the Constraint}\label{direct_min_sec}
Someone may doubt the necessity of the constraint in problem (\ref{non_convex_problem}) and they wonder the performance of the classical accelerated gradient method to minimize problem (\ref{non_convex_problem0}) directly. In this case, the classical accelerated gradient method \cite{Nesterov1988,Nesterov1983,Tseng-2008} becomes
\begin{eqnarray}
&&\V^{k}=(1-\theta_k)\U^{k}+\theta_k\Z^{k},\label{V_step2}\\
&&\Z^{k+1}=\Z^{k}-\eta\nabla g(\V^k),\label{Z_step2}\\
&&\U^{k+1}=(1-\theta_k)\U^{k}+\theta_k\Z^{k+1},\label{U_step2}
\end{eqnarray}
and it is equivalent to
\begin{eqnarray}
&&\V^k=\U^k+\beta_k(\U^k-\U^{k-1}),\label{V_step3}\\
&&\U^{k+1}=\V^k-\eta\nabla g(\V^{k}).\label{U_step3}
\end{eqnarray}
where $\beta_k$ is defined in Lemma \ref{global_lemma}. Another choice is a constant of $\beta<1$. Theorem \ref{direct_min_theorem} establishes the convergence rate for the above two recursions. We leave the proof in Appendix D.
\begin{theorem}\label{direct_min_theorem}
Assume that $\U^*\in\mathcal{X}^*$ and $\V^k\in\mathcal{R}^{n\times r}$ satisfy $\|\V^k-P_{\mathcal{X}^*}(\V^k)\|_F\leq \min\left\{0.01\sigma_r(\U^*), \frac{\mu\sigma_r^2(\U^*)}{6L\|\U^*\|_2}\right\}$. Let $\eta$ be the one in Lemma \ref{global_lemma}. Then, we can have
\begin{eqnarray}
&&g(\U^{k+1})+\nu\|\U^{k+1}-\U^k\|_F^2-g(\U^*)\notag\\
&\leq& \frac{1}{1+\frac{\gamma}{\frac{5}{\eta^2\mu\sigma_r^2(\U^*)}+\nu}}\left[g(\U^k)+\nu\|\U^k-\U^{k-1}\|_F^2-g(\U^*)\right].\notag
\end{eqnarray}
where $\gamma=\frac{1-\beta_{\max}^2}{4\eta}-\frac{\beta_{\max}\hat L}{2}-\frac{\hat L}{4}>0$ and $\nu=\frac{1+\beta_{\max}^2}{4\eta}-\frac{\hat L}{4}>0$.
\end{theorem}

Consider the case that $\beta_k$ is a constant. Then, we know that all of the constants $\gamma,\nu,\hat L$ and $\frac{1}{\eta}$ are of the order $O\left(L\|\U^*\|_2^2+\|\nabla f(\X^*)\|_2\right)$. Thus, the convergence rate of recursion (\ref{V_step3})-(\ref{U_step3}) is in the form of
\begin{eqnarray}
\left(1-\frac{\mu\sigma_r^2(\U^*)}{L\|\U^*\|_2^2+\|\nabla f(\X^*)\|_2}\right)^N,\notag
\end{eqnarray}
which is the same as that of the gradient descent method in (\ref{gd_rate}). Thus, although the convergence of the classical accelerated gradient method for problem (\ref{non_convex_problem0}) can be proved, it is not easy to build the acceleration upon the gradient descent. As a comparison, Algorithm \ref{agd_alg} has a theoretical better dependence on the condition number of $\frac{L}{\mu}$. Thus, the reformulation of problem (\ref{non_convex_problem0}) to a constrained one is necessary to prove acceleration.

\section{The Asymmetric Case}\label{asy_sec}
In this section, we consider the asymmetric case of problem (\ref{convex_problem}):
\begin{eqnarray}
\min_{\widetilde\X\in\mathcal{R}^{n\times m}} f(\widetilde\X),\label{asy_problem}
\end{eqnarray}
where there exists a minimizer $\widetilde\X^*$ of rank-$r$. We follow \cite{Park-2016-siam} to assume $\nabla f(\widetilde\X^*)=0$. In the asymmetric case, we can factorize $\widetilde\X=\widetilde\U\widetilde\V^T$ and reformulate problem (\ref{asy_problem}) as a similar problem to (\ref{non_convex_problem0}). Moreover, we follow \citep{Park-2016-siam,wang-2016} to regularize the objective and force the solution pair $(\widetilde\U,\widetilde\V)$ to be balanced. Otherwise, the problem may be ill-conditioned since $\left(\frac{1}{\delta}\widetilde\U\right)\hspace*{-0.05cm}(\delta\widetilde\V)$ is also a factorization of $\widetilde\U\widetilde\V^T$ for any large $\delta$ \cite{Park-2016-siam}. Specifically, we consider the following problem
\begin{eqnarray}
\min_{\widetilde\U\in\mathcal{R}^{n\times r},\widetilde\V\in\mathcal{R}^{m\times r}} f(\widetilde\U\widetilde\V^T)+\frac{\mu}{8}\|\widetilde\U^T\widetilde\U-\widetilde\V^T\widetilde\V\|_F^2.\label{asy_problem2}
\end{eqnarray}
Let $\widetilde\X^*=\A\Sigma\B^T$ be its SVD. Then, $(\widetilde\U^*=\A\sqrt{\Sigma},\widetilde\V^*=\B\sqrt{\Sigma})$ is a minimizer of problem (\ref{asy_problem2}). Define a stacked matrix $\U=\left(
  \begin{array}{c}
     \widetilde\U \\
     \widetilde\V
  \end{array}
\right)$ and let $\X=\U\U^T=\left(
  \begin{array}{cc}
     \widetilde\U\widetilde\U^T & \widetilde\U\widetilde\V^T\\
     \widetilde\V\widetilde\U^T & \widetilde\V\widetilde\V^T
  \end{array}
\right)$. Then we can write the objective in (\ref{asy_problem2}) in the form of $\hat f(\X)$, defined as $\hat f(\X)=f(\widetilde\U\widetilde\V^T)+\frac{\mu}{8}\|\widetilde\U\widetilde\U^T\|_F^2+\frac{\mu}{8}\|\widetilde\V\widetilde\V^T\|_F^2-\frac{\mu}{4}\|\widetilde\U\widetilde\V^T\|_F^2$. Since $f$ is restricted $\mu$-strongly convex, we can easily check that $\hat f(\X)$ is restricted $\frac{\mu}{4}$-strongly convex. On the other hand, we know that $\hat f(\X)$ is restricted $\left(L+\frac{\mu}{2}\right)$-smooth. Applying the conclusions on the symmetric case to $\hat f(\X)$, we can apply Algorithm \ref{agd_alg} to the asymmetric case. From Theorem \ref{conver_ratem}, we can get the convergence rate. Moreover, since $\sigma_i(\X^*)=2\sigma_i(\widetilde\X^*)$,
\begin{eqnarray}
\begin{aligned}
\nabla f(\X^*)=&\left(
  \begin{array}{cc}
     0 & \nabla f(\widetilde\X^*)\\
     \nabla f(\widetilde\X^*)^T & 0
  \end{array}
\right)+\frac{\mu}{4}\left(
  \begin{array}{c}
     \widetilde\U^* \\
     -\widetilde\V^*
  \end{array}
\right)\left(\widetilde{\U^*}^T,-\widetilde{\V^*}^T\right)\\
=&\frac{\mu}{4}\left(
  \begin{array}{c}
     \widetilde\U^* \\
     -\widetilde\V^*
  \end{array}
\right)\left(\widetilde{\U^*}^T,-\widetilde{\V^*}^T\right)\notag
\end{aligned}
\end{eqnarray}
and $\|\nabla f(\X^*)\|_2=\frac{\mu}{4}\|\X^*\|_2$, where $\X^*=\left(
  \begin{array}{cc}
     \widetilde\U^*\widetilde{\U^*}^T & \widetilde\U^*\widetilde{\V^*}^T\\
     \widetilde\V^*\widetilde{\U^*}^T & \widetilde\V^*\widetilde{\V^*}^T
  \end{array}
\right)$ and we use $\nabla f(\widetilde\X^*)=0$, we can simplify the worst case convergence rate to
$\left(1-\frac{\sigma_r(\widetilde\X^*)}{\|\widetilde\X^*\|_2}\sqrt{\frac{\mu}{(m+n)rL}}\right)^N$. As a comparison, the rate of the gradient descent is $\left(1-\frac{\sigma_r(\widetilde\X^*)}{\|\widetilde\X^*\|_2}\frac{\mu}{L}\right)^N$ \cite{Park-2016-siam}.

In the asymmetric case, both $\widetilde\U^*$ and $\widetilde\V^*$ are of full rank. Otherwise, $\rank(\widetilde\X^*)<r$. Thus, we can select the index set $ S^1$ from $\widetilde\U^0$ and select $ S^2$ from $\widetilde\V^0$ with the guarantee of $\sigma_r(\widetilde\U^0_{ S^1})\geq \frac{\sigma_r(\widetilde\U^0)}{\sqrt{2r(n-r+1)}}$ and $\sigma_r(\widetilde\V^0_{ S^2})\geq \frac{\sigma_r(\widetilde\V^0)}{\sqrt{2r(m-r+1)}}$.

\section{Experiments}\label{mc_exp}

In this section, we test the efficiency of the proposed Accelerated Gradient Descent (AGD) method on Matrix Completion, One Bit Matrix Completion and Matrix Regression.

\subsection{Matrix Completion}

In matrix completion \citep{Rohde-2011,Koltchinsii-2011,Negahban-2012}, the goal is to recover the low rank matrix $\X^*$ based on a set of randomly observed entries $\mathbf O$ from $\X^*$. The traditional matrix completion problem is to solve the following model:
\begin{eqnarray}
\min_{\X} \frac{1}{2}\sum_{(i,j)\in\mathbf O}(\X_{i,j}-\X_{i,j}^*)^2,\quad s.t.\quad \rank(\X)\leq r.\notag
\end{eqnarray}
We consider the asymmetric case and solve the following model:
\begin{eqnarray}
\min_{\widetilde\U\in\mathcal{R}^{n\times r},\widetilde\V\in\mathcal{R}^{m\times r}} \frac{1}{2}\sum_{(i,j)\in\mathbf O}((\widetilde\U\widetilde\V^T)_{i,j}-\X_{i,j}^*)^2+\frac{1}{200}\|\widetilde\U^T\widetilde\U-\widetilde\V^T\widetilde\V\|_F^2.\notag
\end{eqnarray}

\begin{figure}
\centering
\begin{tabular}{@{\extracolsep{0.001em}}c@{\extracolsep{0.001em}}c@{\extracolsep{0.001em}}c@{\extracolsep{0.001em}}c}
\includegraphics[width=0.33\textwidth,keepaspectratio]{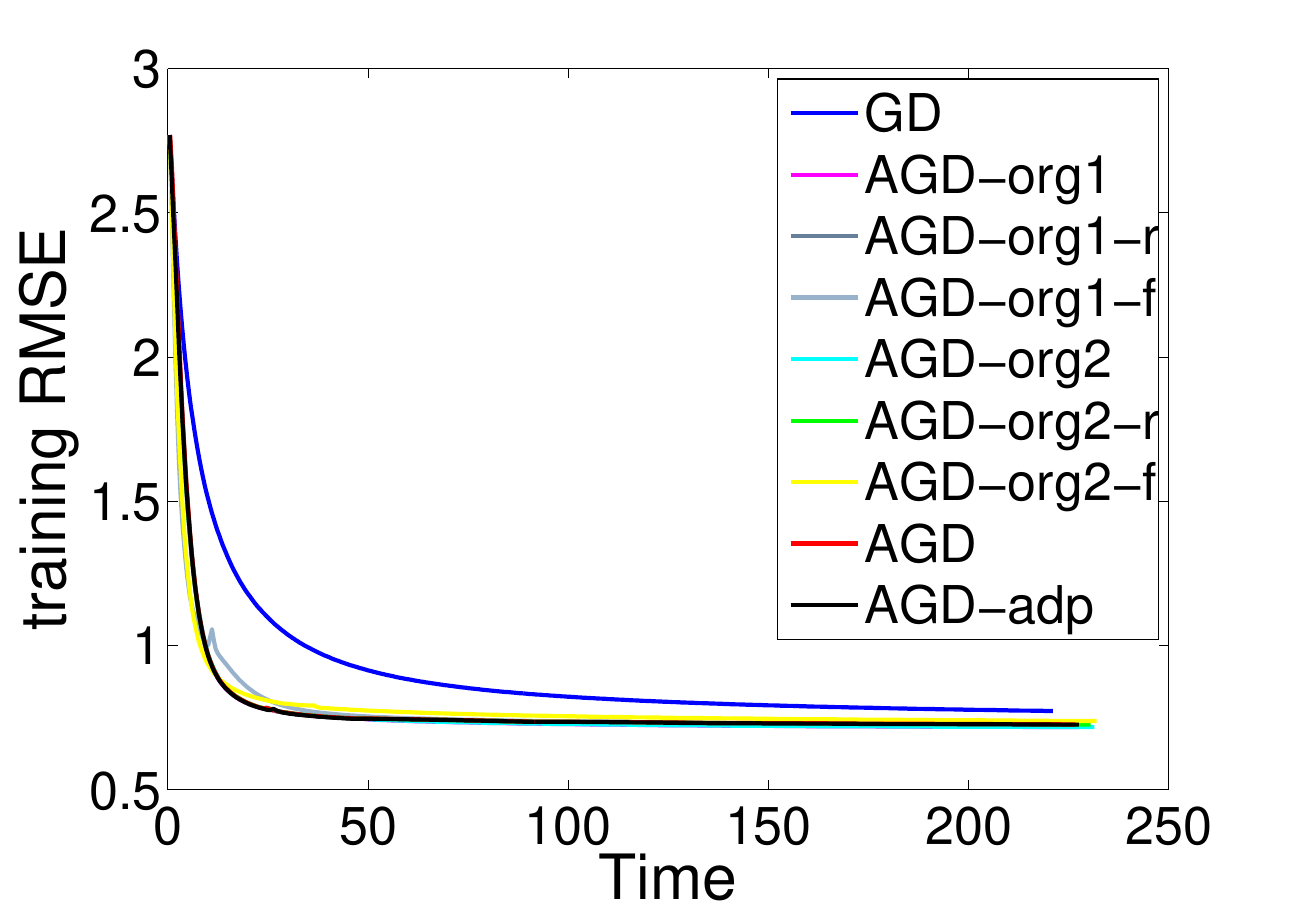}
&\includegraphics[width=0.33\textwidth,keepaspectratio]{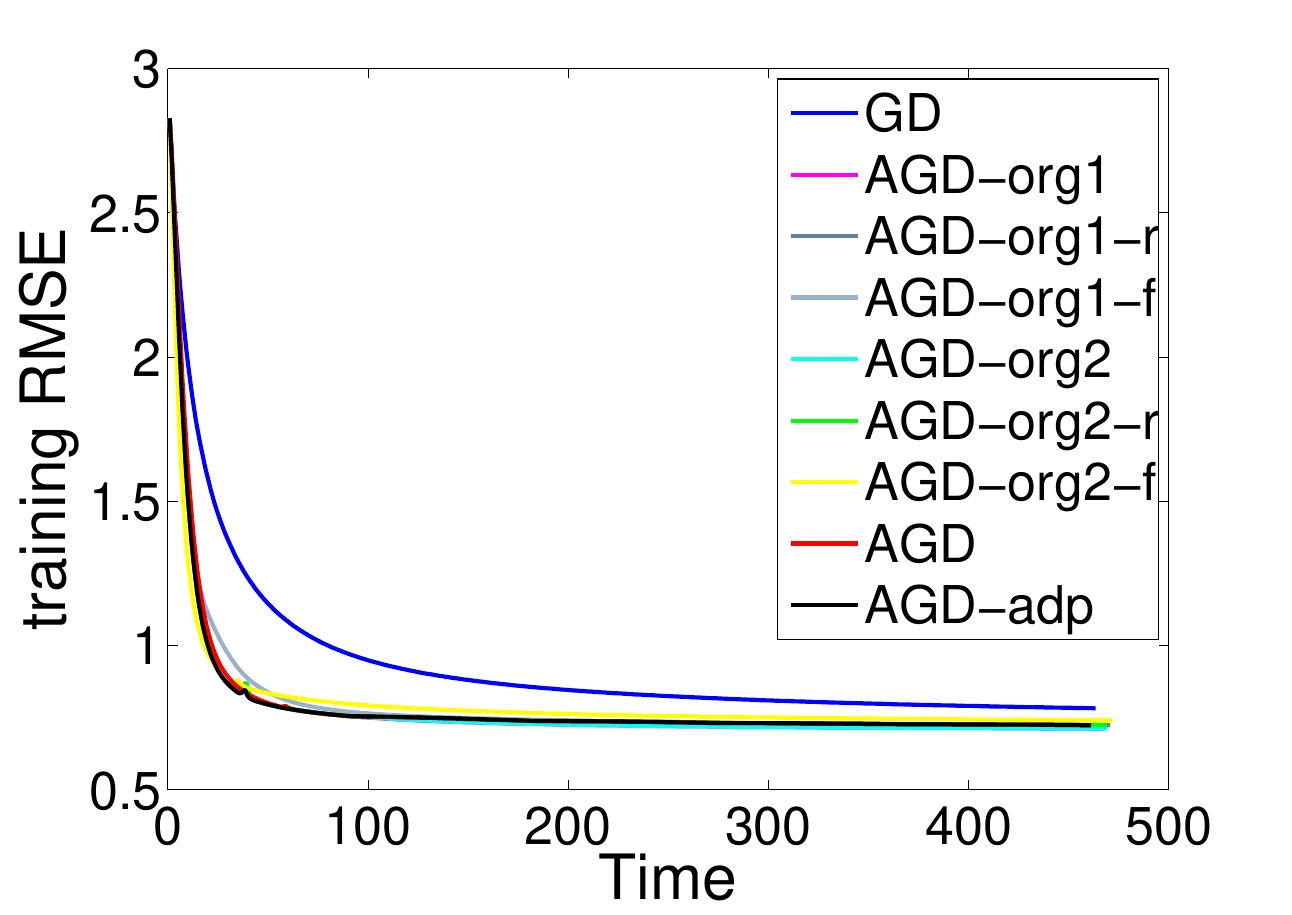}
&\includegraphics[width=0.33\textwidth,keepaspectratio]{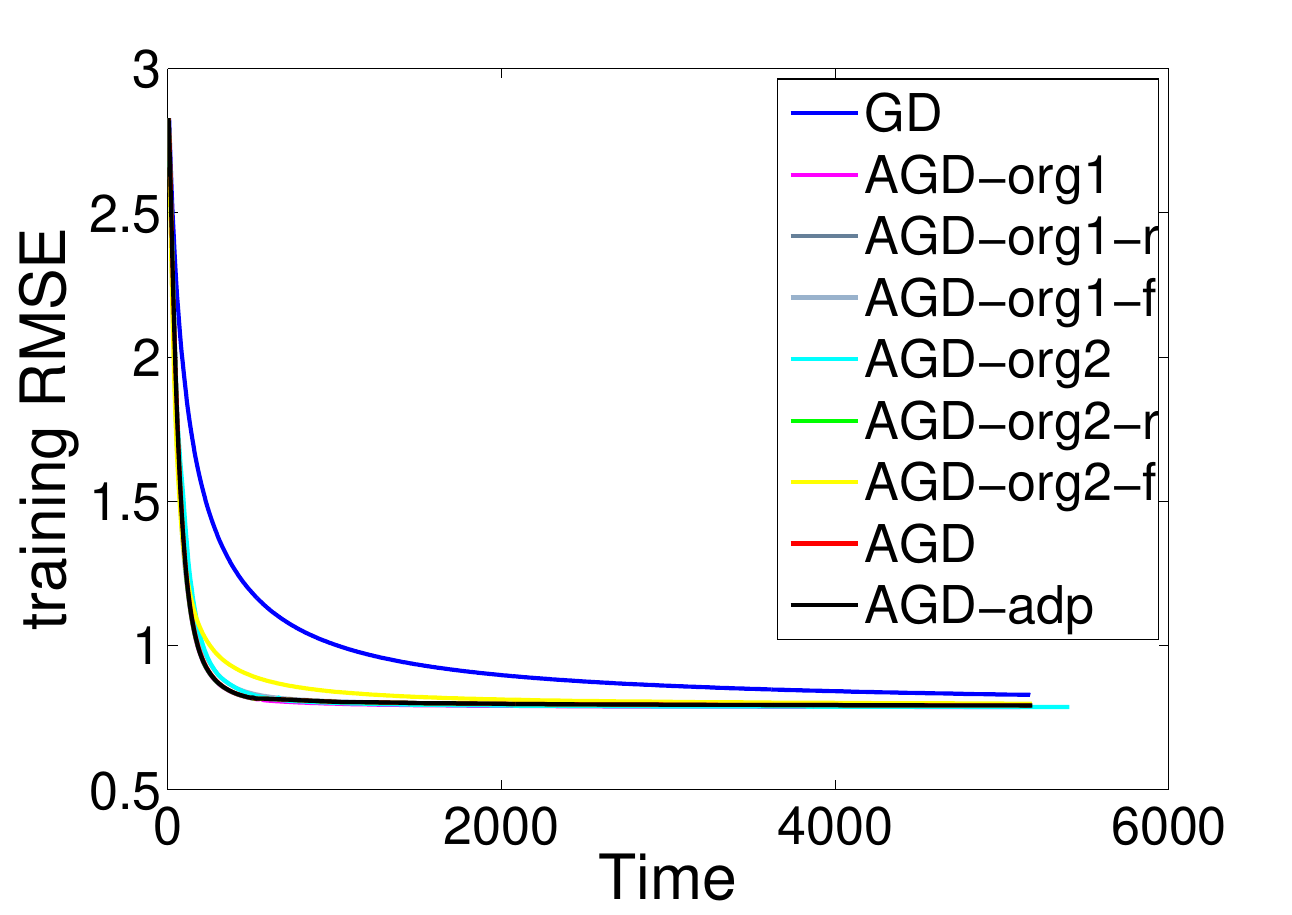}\\
\includegraphics[width=0.33\textwidth,keepaspectratio]{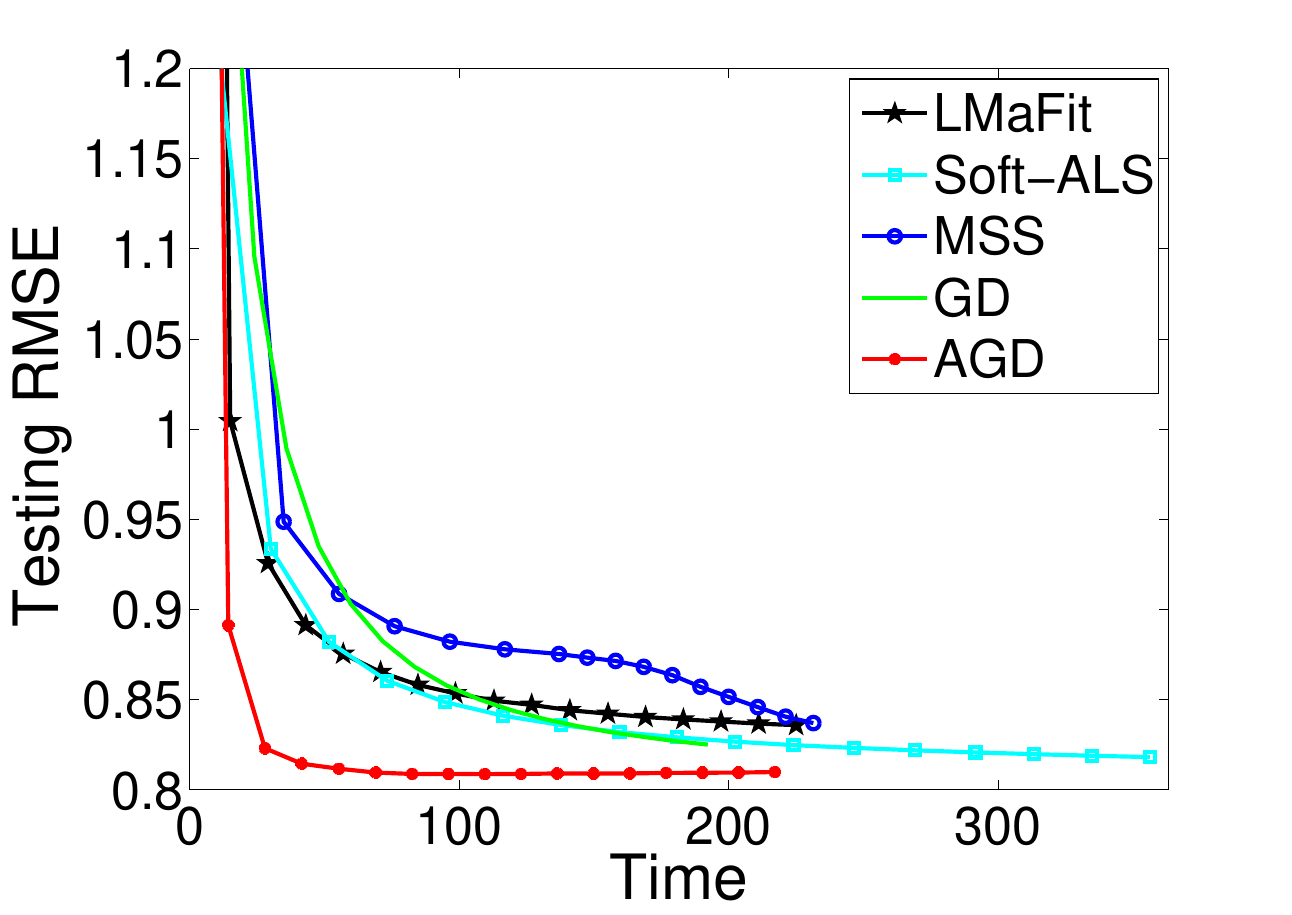}
&\includegraphics[width=0.33\textwidth,keepaspectratio]{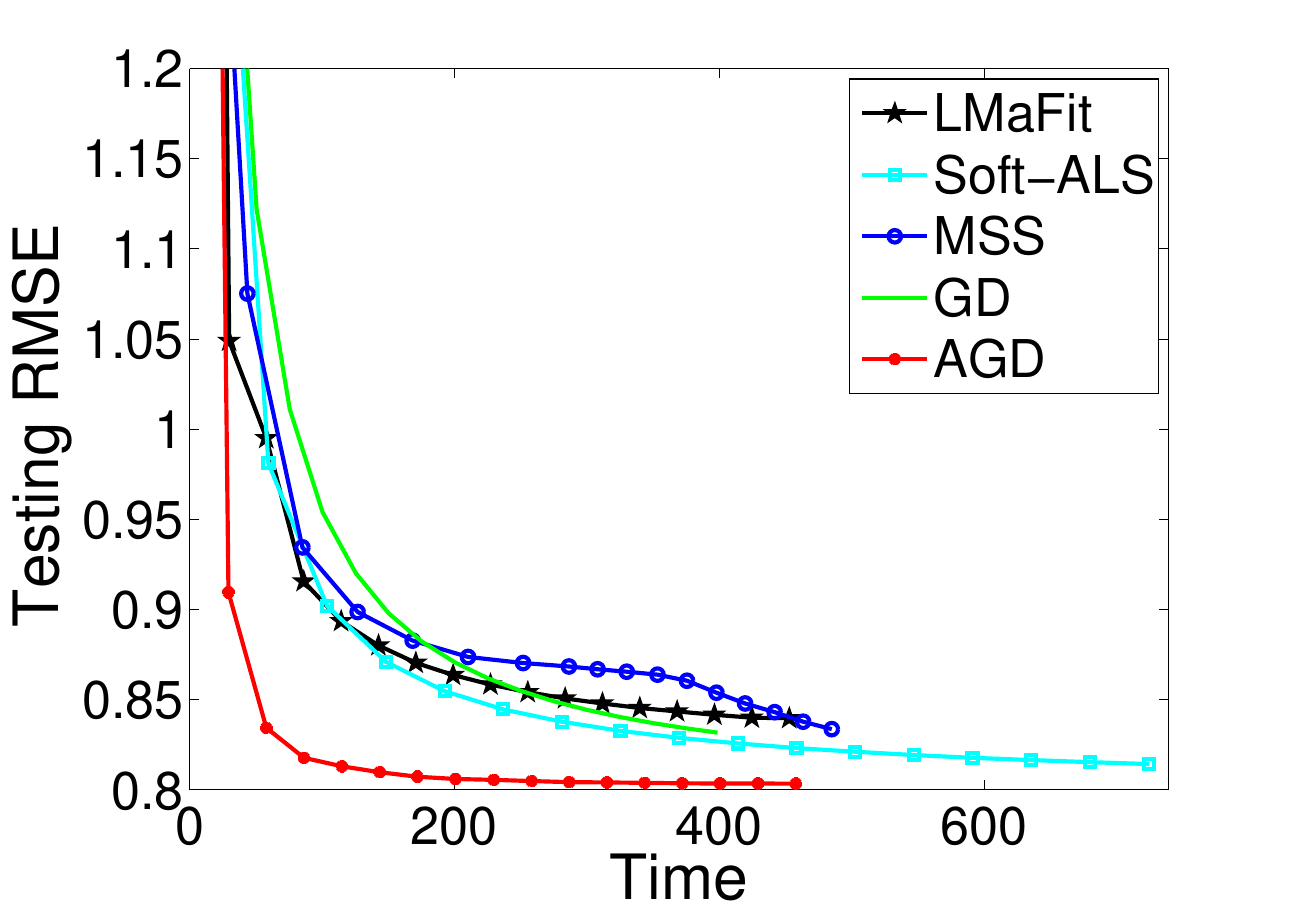}
&\includegraphics[width=0.33\textwidth,keepaspectratio]{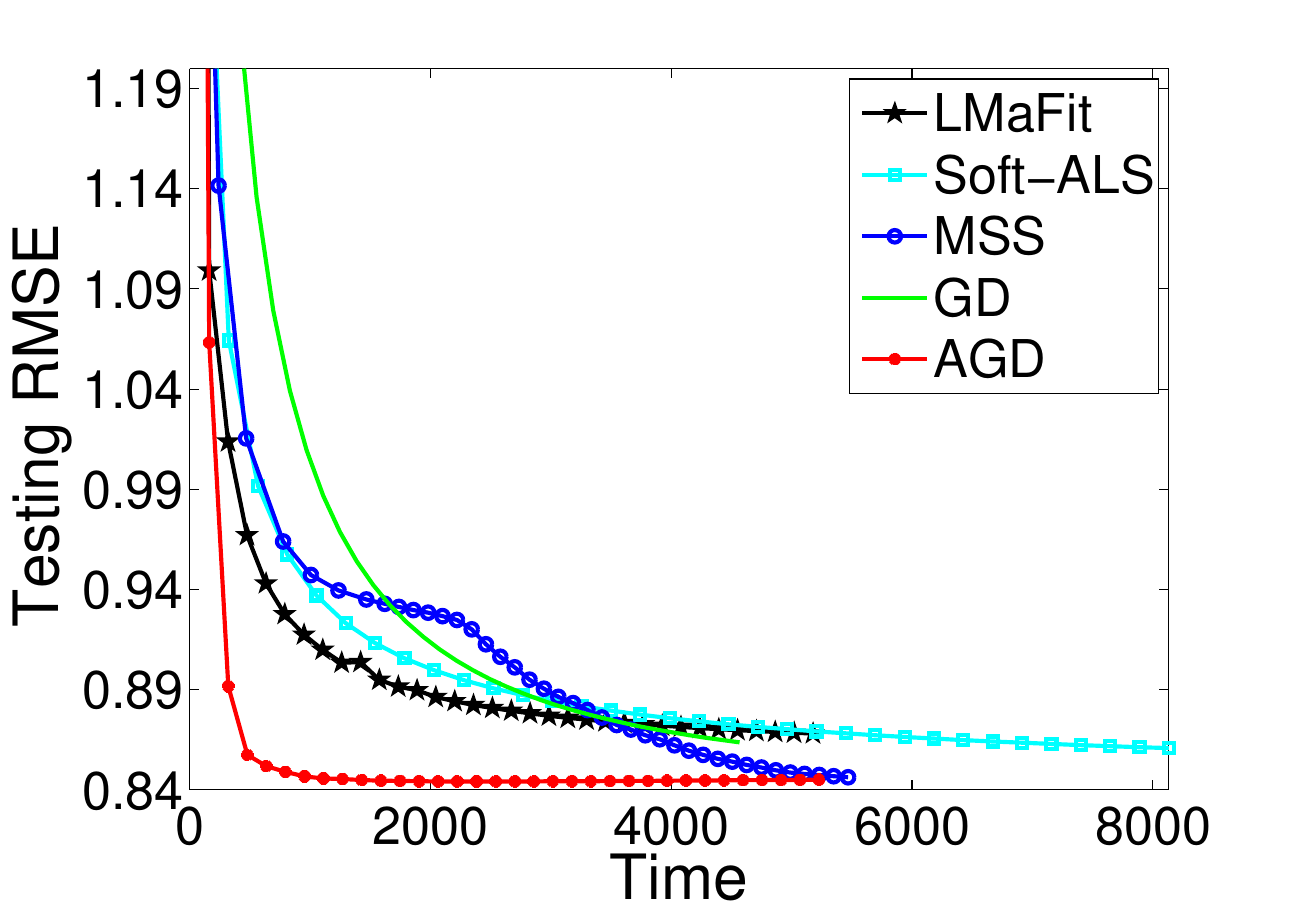}\\
(a) MovieLens-10M & (b) MovieLens-20M & (c) Netflix
\end{tabular}
\caption{Top: Compare the training RMSE of GD, AGD, AGD-adp and several variants of the original AGD. Bottom: Compare the testing RMSE of GD, AGD, LMaFit, Soft-ALS and MSS.}\label{figC3}
\end{figure}

We set $r=10$ and test the algorithms on the Movielen-10M, Movielen-20M and Netflix data sets. The corresponding observed matrices are of size $69878\times 10677$ with $o\%=1.34\%$, $138493\times 26744$ with $o\%=0.54\%$ and $480189\times 17770$ with $o\%=1.18\%$, respectively, where $o\%$ means the percentage of the observed entries. We compare AGD and AGD-adp (AGD with adaptive index sets selection) with GD and several variants of the original AGD: \newline
1. AGD-original1: The classical AGD with recursions of (\ref{V_step3})-(\ref{U_step3}).\newline
2. AGD-original1-r: AGD-original1 with restart.\newline
3. AGD-original1-f: AGD-original1 with fixed $\beta_k$ of $\frac{\sqrt{L}-\sqrt{\mu}}{\sqrt{L}+\sqrt{\mu}}$.\newline
4. AGD-original2: The classical AGD with recursions of (\ref{V_step2})-(\ref{U_step2}).\newline
5. AGD-original2-r: AGD-original2 with restart.\newline
6. AGD-original2-f: AGD-original2 with fixed $\theta$.

Let $\X_{\mathbf O}$ be the observed data and $\A\bm{\Sigma}\B^T$ be its SVD. We initialize $\widetilde\U=\A_{:,1:r}\sqrt{\bm{\Sigma}_{1:r,1:r}}$ and $\widetilde\V=\B_{:,1:r}\sqrt{\bm{\Sigma}_{1:r,1:r}}$ for all the compared methods. Since $\X_{\mathbf O}$ is sparse, it is efficient to find the top $r$ singular values and the corresponding singular vectors for large scale matrices \citep{lansvd}. We tune the best step sizes of $\eta=5\times 10^{-5},4\times 10^{-5}$ and $1\times 10^{-5}$ for all the compared methods on the three data sets, respectively. For AGD, we set $\epsilon=10^{-10}$, $ S^1=\{1:r\}$ and $ S^2=\{r+1:2r\}$ for simplicity. We set $K = 100$ for AGD, AGD-adp and the original AGD with restart. We run the compared methods 500 iterations for the Movielen-10M and Movielen-20M data sets and 1000 iterations for the Netflix data set.

The top part of Figure \ref{figC3} plots the curves of the training RMSE v.s. time (seconds). We can see that AGD is faster than GD. The performances of AGD, AGD-adp and the original AGD are similar. In fact, in AGD-adp, we observe that the index sets do not change during the iterations. Thus, the condition of $\sigma_r(\U^{t,K+1}_{ S'})\geq\epsilon$ $\forall t$ in Theorem \ref{global_convergencem} holds. The original AGD performs almost equally fast as our modified AGD in practice. However, it has an inferior convergence rate theoretically. The bottom part of Figure \ref{figC3} plots the curves of the testing RMSE v.s. time. Besides GD, we also compare AGD with LMaFit \citep{Wen-2012}, Soft-ALS \citep{Hastie-2015} and MSS \citep{Xu-2016}. They all solve a factorization based nonconvex model. From Figure \ref{figC3} we can see that AGD achieves the lowest testing RMSE with the fastest speed.

\subsection{One Bit Matrix Completion}

In one bit matrix completion \citep{Davenport-2014}, the sign of a random subset from the unknown low rank matrix $\X^*$ is observed, instead of observing the actual entries. Given a probability density function $f$, e.g., the logistic function $f(\x)=\frac{e^x}{1+e^x}$, we observe the sign of $\x$ as $+1$ with probability $f(\x)$ and observe the sign as $-1$ with probability $1-f(\x)$. The training objective is to minimize the negative log-likelihood:
\begin{eqnarray}
\min_{\X} -\sum_{(i,j)\in\mathbf O}\{\mathbf 1_{\Y_{i,j}=1}\mbox{log}(f(\X_{i,j}))+\mathbf 1_{\Y_{i,j}=-1}\mbox{log}(1-f(\X_{i,j}))\},s.t. \rank(\X)\leq r.\notag
\end{eqnarray}

\begin{figure}
\centering
\begin{tabular}{@{\extracolsep{0.001em}}c@{\extracolsep{0.001em}}c@{\extracolsep{0.001em}}c@{\extracolsep{0.001em}}c}
\includegraphics[width=0.33\textwidth,keepaspectratio]{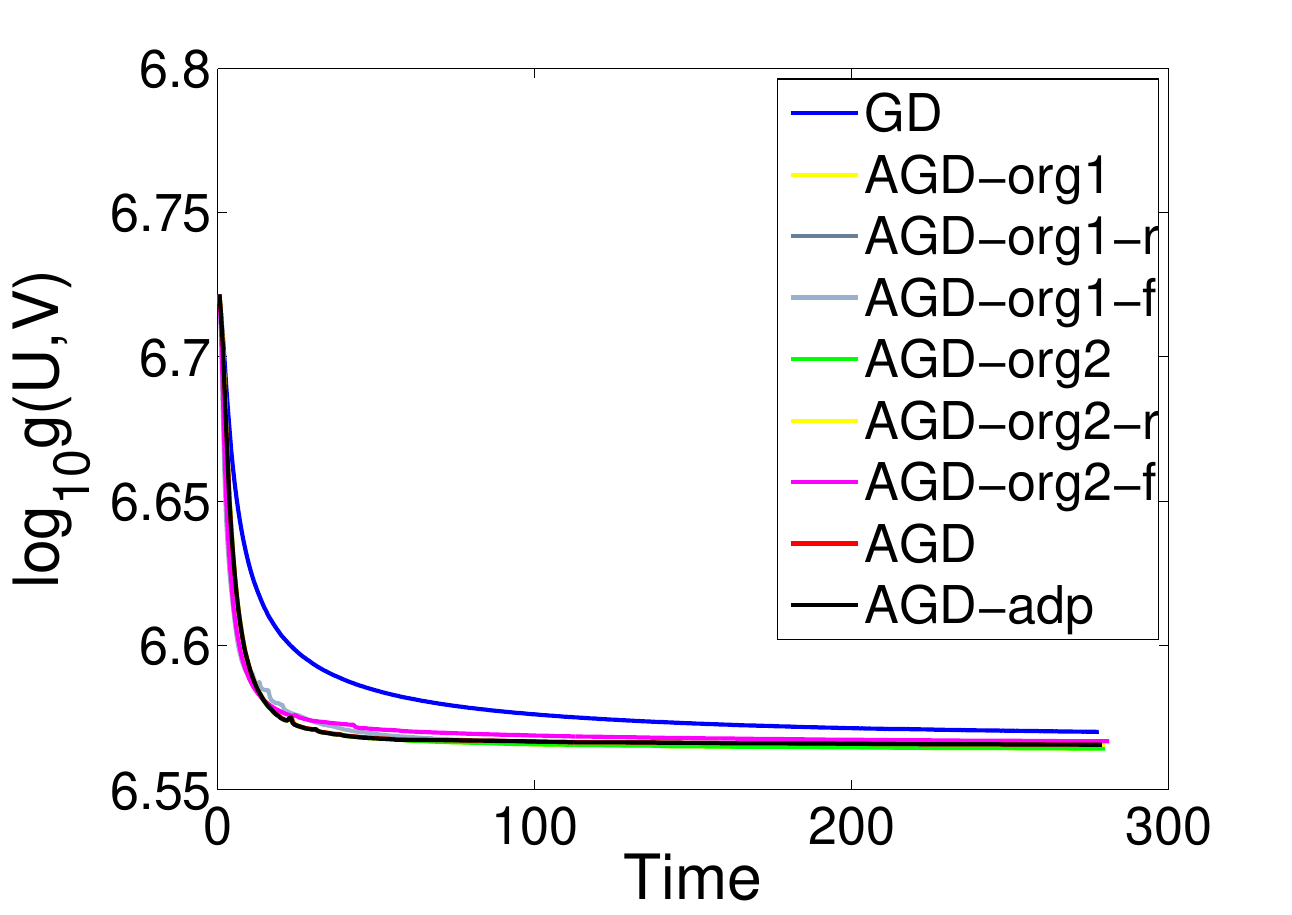}
&\includegraphics[width=0.33\textwidth,keepaspectratio]{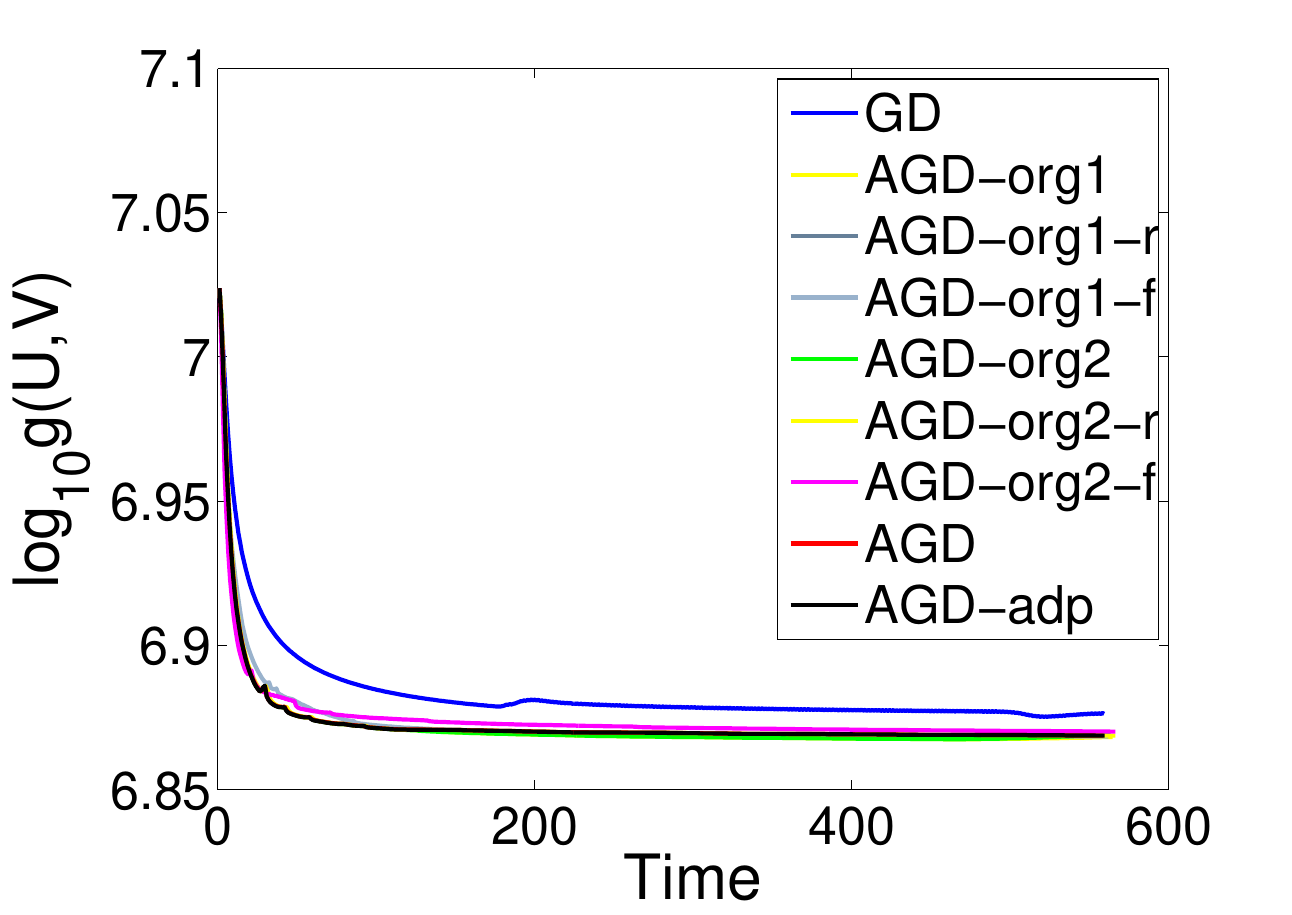}
&\includegraphics[width=0.33\textwidth,keepaspectratio]{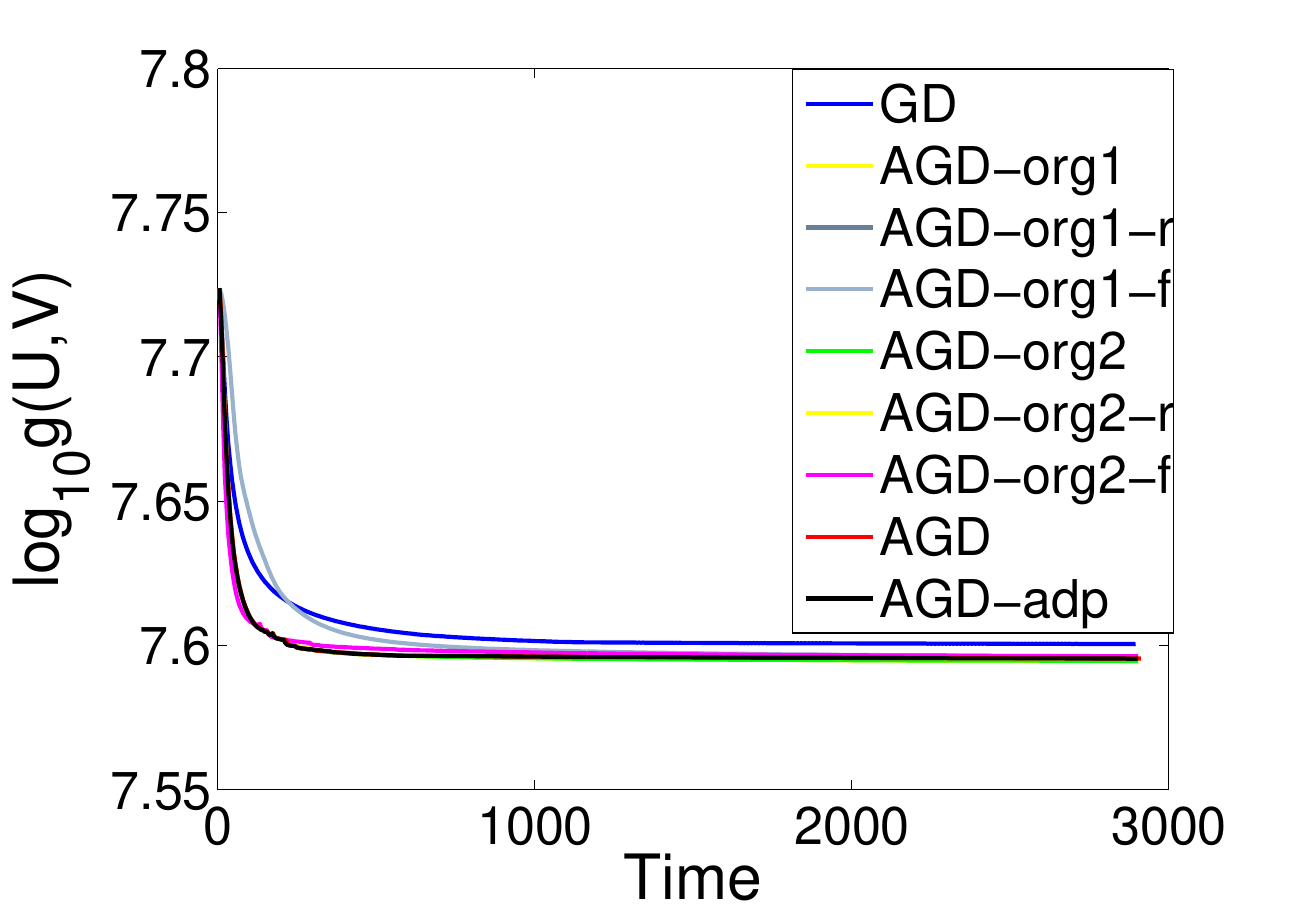}\\
(a) MovieLens-10M & (b) MovieLens-20M & (c) Netflix
\end{tabular}
\caption{Compare AGD and AGD-adp with GD and several variants of the original AGD on the One Bit Matrix Completion problem.}\label{figB2}
\end{figure}

In this section, we solve the following model:
\begin{eqnarray}
\min_{\widetilde\U,\widetilde\V}&&-\sum_{(i,j)\in\mathbf O}\left\{\mathbf 1_{\Y_{i,j}=1}\mbox{log}(f((\widetilde\U\widetilde\V^T)_{i,j}))+\mathbf 1_{\Y_{i,j}=-1}\mbox{log}(1-f((\widetilde\U\widetilde\V^T)_{i,j}))\right\}\notag\\
&&+\frac{1}{200}\|\widetilde\U^T\widetilde\U-\widetilde\V^T\widetilde\V\|_F^2.\notag
\end{eqnarray}
We use the data sets of Movielen-10M, Movielen-20M and Netflix. We set $\Y_{i,j}=1$ if the $(i,j)$-th observation is larger than the average of all observations and $\Y_{i,j}=-1$, otherwise. We set $r=5$ and $\eta=0.001$, $0.001, 0.0005$ for all the compared methods on the three data sets. The other experimental setting is the same as Matrix Completion. We run all the methods for 500 iterations. Figure \ref{figB2} plots the curves of the objective value v.s. time (seconds) and we can see that AGD is also faster than GD. The performances of AGD, AGD-adp and the original AGD are nearly the same.

\subsection{Matrix Regression}

\begin{figure}
\centering
\begin{tabular}{@{\extracolsep{0.001em}}c@{\extracolsep{0.001em}}c@{\extracolsep{0.001em}}c@{\extracolsep{0.001em}}c}
\includegraphics[width=0.33\textwidth,keepaspectratio]{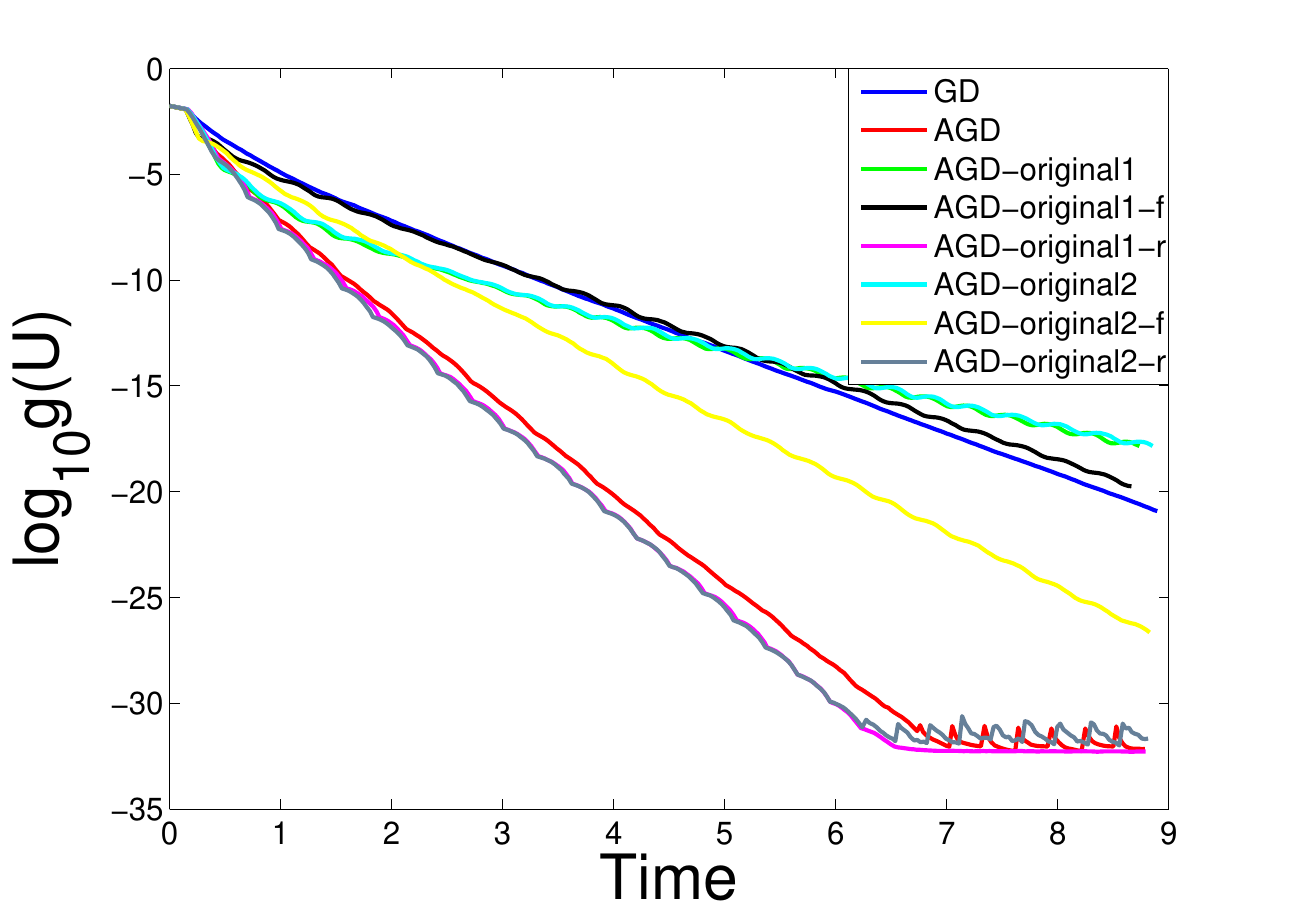}
&\includegraphics[width=0.33\textwidth,keepaspectratio]{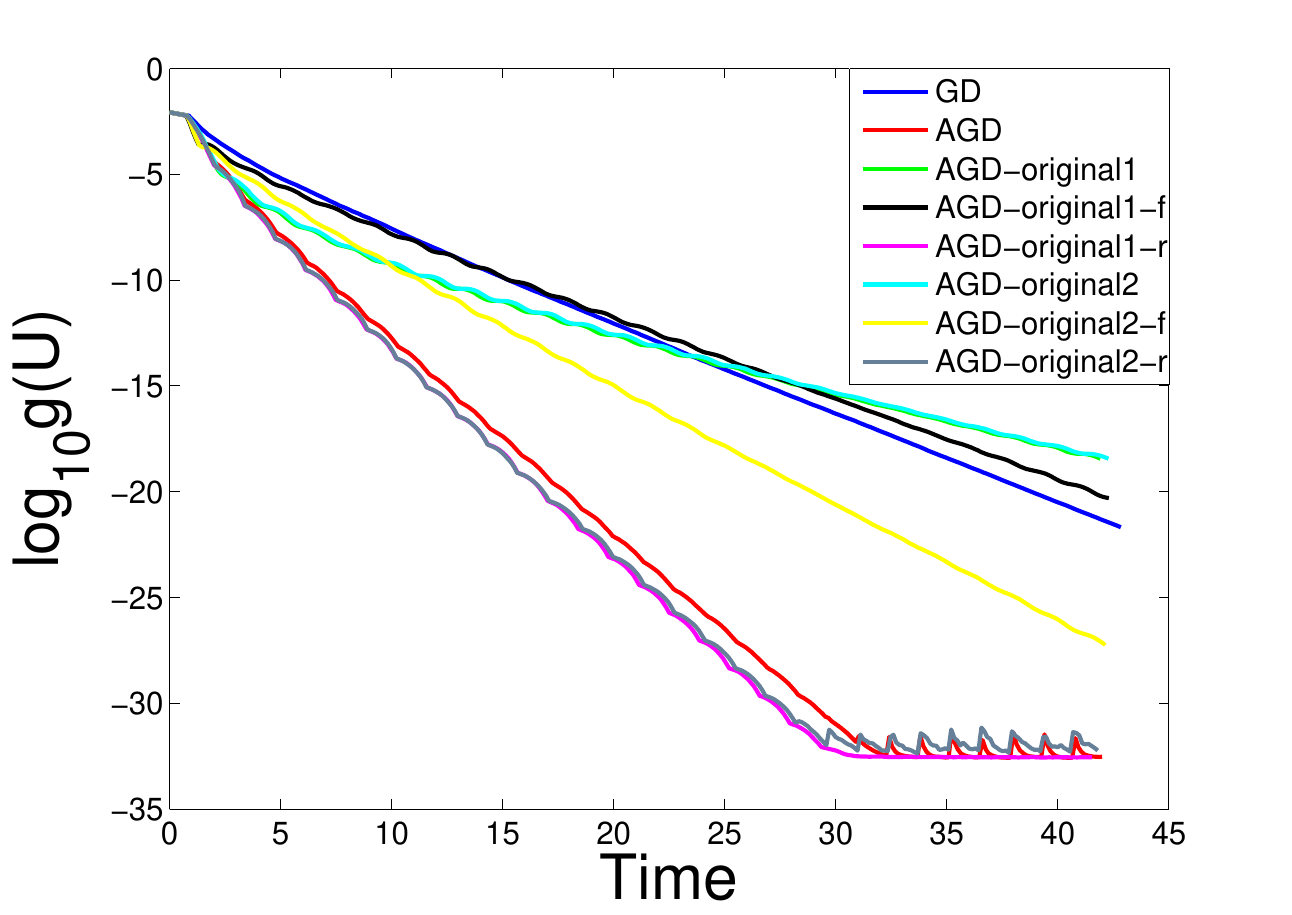}
&\includegraphics[width=0.33\textwidth,keepaspectratio]{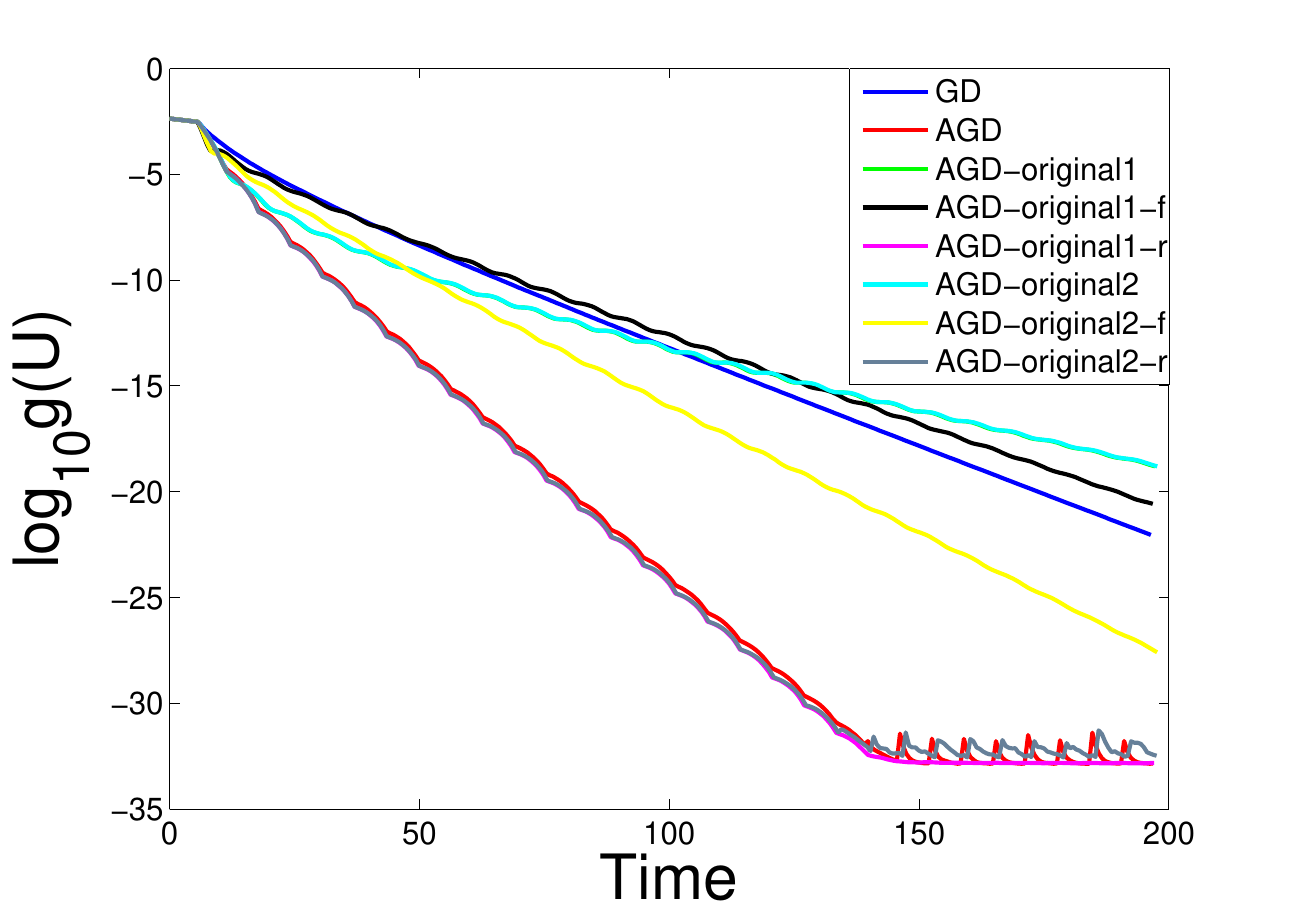}\\
(a) $n$=512 & (b) $n$=1024 & (c) $n$=2048
\end{tabular}
\caption{Compare AGD with GD and several variants of the original AGD on the Matrix Regression problem.}\label{figR}
\end{figure}

In matrix regression \citep{Recht-2010,Negahban-2011}, the goal is to estimate the unknown low rank matrix $\X^*$ from a set of measurements $\y=\A(\X^*)+\varepsilon$, where $\A$ is a linear operator and $\varepsilon$ is the noise. A reasonable estimation of $\X^*$ is to solve the following rank constrained problem:
\begin{eqnarray}
\min_{\X}f(\X)=\frac{1}{2}\|\A(\X)-\y\|_F^2,\quad s.t. \quad \rank(\X)\leq r.\notag
\end{eqnarray}
We consider the symmetric case of $\X$ and solve the following nonconvex model:
\begin{eqnarray}
\min_{\U\in\mathcal{R}^{n\times r}} f(\U)=\frac{1}{2}\|\A(\U\U^T)-\y\|_F^2.\notag
\end{eqnarray}
We follow \citep{Srinadh-2016-colt} to use the permuted and sub-sampled noiselets \citep{Waters-noisylet} for the linear operator $\A$ and $\U^*$ is generated from the normal Gaussian distribution without noise. We set $r=10$ and test different $n$ with $n$=512, 1024 and 2048. We fix the number of measurements to $4nr$ and follow \citep{Srinadh-2016-colt} to use the initializer from the eigenvalue decomposition of $\frac{\X^0+(\X^0)^T}{2}$ for all the compared methods, where $\X^0=\mbox{Project}_{+}\left(\frac{-\nabla f(0)}{\|\nabla f(0)-\nabla f(11^T)\|_F}\right)$. We set $\eta=5,10$ and $20$ for all the compared methods for $n=512,1024$ and $2048$, respectively. In AGD, we set $\epsilon=10^{-10}$, $K=10$, $ S^1=\{1:r\}$ and $ S^2=\{r+1:2r\}$. Figure \ref{figR} plots the curves of the objective value v.s. time (seconds). We run all the compared methods for 300 iterations. We can see that AGD and the original AGD with restart perform almost equally fast. AGD runs faster than GD and the original AGD without restart.

\subsection{Verifying (\ref{cont15}) in Practice}\label{section_exp_4}

In this section, we verify that the conditions of $\|\hat\U_S^{t,k}-\U_S^*\|_F\leq c\sqrt{\frac{r}{n}}\|\hat\U^{t,k}-\U^*\|_F$ and $\|\hat\V_S^{t,k}-\U_S^*\|_F\leq c\sqrt{\frac{r}{n}}\|\hat\V^{t,k}-\U^*\|_F$ in (\ref{cont15}) hold in our experiments, where $\hat\U^{t,k}=\U^{t,k}\R$ with $\R=\argmin_{\R\mathcal{R}^{r\times r},\R\R^T=\I}\|\U^{t,k}\R-\U^*\|_F^2$ and $\hat\V^{t,k}$ is similar. We use the final output $\U^{T,K+1}$ as $\U^*$. Table \ref{table-cc} lists the results. We can see that $\frac{\|\hat\U^{t,k}_{S}-\U^{*}_{S}\|_F}{\|\hat\U^{t,k}-\U^{*}\|_F}$ and $\frac{\|\hat\V^{t,k}_{S}-\U^{*}_{S}\|_F}{\|\hat\V^{t,k}-\U^{*}\|_F}$ have the same order as $\sqrt{\frac{r}{n}}$.

\begin{table}
\caption{Testing the order of $\frac{\|\hat\U^{t,k}_{S}-\U^{*}_{S}\|_F}{\|\hat\U^{t,k}-\U^{*}\|_F}$ and $\frac{\|\hat\V^{t,k}_{S}-\U^{*}_{S}\|_F}{\|\hat\V^{t,k}-\U^{*}\|_F}$.}\label{table-cc}
\begin{center}
\begin{tabular}{c|c||c|c|c|c|c|c|c}\hline
\multirow{2}{*}{Problem} & \multirow{2}{*}{Data} & \multicolumn{3}{c|}{$\frac{\|\hat\U^{t,k}_{S}-\U^{*}_{S}\|_F}{\|\hat\U^{t,k}-\U^{*}\|_F}$} & \multicolumn{3}{c|}{$\frac{\|\hat\V^{t,k}_{S}-\U^{*}_{S}\|_F}{\|\hat\V^{t,k}-\U^{*}\|_F}$} & \multirow{2}{*}{$\sqrt{\frac{r}{n}}$}\\
\cline{3-8}
& & max & average & min & max & average & min &
\\\hline\hline
\multirow{3}{*}{MR}   & 512  & 0.1536  & 0.1521  & 0.1505 & 0.1536  & 0.1521  & 0.1505  & 0.1398\\
                      & 1024 & 0.0984  & 0.0939  & 0.0894 & 0.0984  & 0.0939  & 0.0894  & 0.0988\\
                      & 2048 & 0.0715  & 0.0681  & 0.0648 & 0.0715  & 0.0681  & 0.0648  & 0.0699\\
                      \hline
\multirow{3}{*}{1bit-MC} & 512  & 0.0344  & 0.0086  & 0.0021 & 0.0344  & 0.0086  & 0.0017 & 0.0079\\
                         & 1024 & 0.0330  & 0.0077  & 0.0022 & 0.0329  & 0.0077  & 0.0020 & 0.0055\\
                         & 2048 & 0.0189  & 0.0103  & 0.0068 & 0.0151  & 0.0103  & 0.0062 & 0.0032\\
                      \hline
\multirow{3}{*}{MC}      & 512  & 0.0664  & 0.0280  & 0.0191 & 0.0664  & 0.0280  & 0.0191 & 0.0111\\
                         & 1024 & 0.0569  & 0.0230  & 0.0151 & 0.0569  & 0.0230  & 0.0139 & 0.0078\\
                         & 2048 & 0.0346  & 0.0191  & 0.0105 & 0.0346  & 0.0190  & 0.0104 & 0.0045\\

 \hline
\end{tabular}
\end{center}
\end{table}

\section{Conclusions}
In this paper, we study the factorization based low rank optimization. A linearly convergent accelerated gradient method with alternating constraint is proposed with the optimal dependence on the condition number of $\sqrt{L/\mu}$ as convex programming. As far as we know, this is the first work with the provable optimal dependence on $\sqrt{L/\mu}$ for this kind of nonconvex problems. Our method can also be applied to the asymmetric factorization.

Many open problems remain in this area. 1. What is the lower complexity bound of problem (\ref{non_convex_problem0})? Is our theoretical upper complexity bound tight on the dependence of $n$ when $L/\mu\geq O(n)$? 2. Can the original AGD have the guaranteed acceleration theoretically? In other words, can we improve Necoara, Nesterov and Glineur \cite{Necoara-2016}'s results without the uniqueness assumption and have a better result than Theorem \ref{direct_min_theorem}?


\section*{Appendix A}\label{sec_appendixA}

\begin{lemma}\label{sdp_lemma}
For problem (\ref{convex_problem}) and its minimizer $\X^*$, we have
\begin{eqnarray}
\nabla f(\X^*)\succeq 0.\notag
\end{eqnarray}
\end{lemma}
\begin{proof}
Introduce the Lagrange function
\begin{eqnarray}
L(\X,\bm{\Lambda})=f(\X)+\<\bm{\Lambda},\X\>.\notag
\end{eqnarray}
Since $\X^*$ is the minimizer of problem (\ref{convex_problem}), we know that there exists $\bm{\Lambda}^*$ such that
\begin{eqnarray}
&&\nabla f(\X^*)+\bm{\Lambda}^*=0,\notag\\
&&\<\bm{\Lambda}^*,\X^*\>=0,\quad \X^*\succeq 0,\quad \bm{\Lambda}^*\preceq 0.\notag
\end{eqnarray}
Thus, we can have the conclusion.
\end{proof}

\begin{lemma}\label{diff_lower_bound}\citep{Tu-2016}
For any $\U\in\mathcal{R}^{n\times r},\V\in\mathcal{R}^{n\times r}$, let $\R=\argmin_{\R\R^T=\I}\|\V\R-\U\|_F^2$ and $\hat\V=\V\R$. Then, we can have
\begin{eqnarray}
\|\V\V^T-\U\U^T\|_F^2\geq (2\sqrt{2}-2)\sigma_r^2(\U)\|\hat\V-\U\|_F^2.\notag
\end{eqnarray}
\end{lemma}

\begin{lemma}\citep{Srinadh-2016-colt}\label{variable_bound}
Assume that $\|\U-\U^*\|_F\leq 0.01\sigma_r(\U^*)$. Then, we can have
\begin{eqnarray}
&0.99\sigma_r(\U^*)\leq \sigma_r(\U)\leq 1.01\sigma_r(\U^*),\notag\\
&0.99\|\U^*\|_2\leq \|\U\|_2\leq 1.01\|\U^*\|_2.\notag
\end{eqnarray}
\end{lemma}

\begin{lemma}\label{matrix_lemma}
For any $\U,\V\in\mathcal{R}^{n\times r}$, we have
\begin{eqnarray}
&\|\U\U^T-\V\V^T\|_F\leq (\|\U\|_2+\|\V\|_2)\|\U-\V\|_F,\label{cont10}\\
&\|\nabla f(\V\V^T)-\nabla f(\U^*{\U^*}^T)\|_2\hspace*{-0.08cm}\leq L(\|\V\|_2+\|\U^*\|_2)\|\V-\U^*\|_F,\label{cont12}\\
&\|\nabla g(\U)\hspace*{-0.08cm}-\hspace*{-0.08cm}\nabla g(\V)\|_F\hspace*{-0.08cm}\leq\hspace*{-0.08cm} \left(L\|\U\|_2(\|\U\|_2\hspace*{-0.08cm}+\hspace*{-0.08cm}\|\V\|_2)\hspace*{-0.08cm}+\hspace*{-0.08cm}\|\nabla f(\V\V^T)\|_2\right)\|\U\hspace*{-0.08cm}-\hspace*{-0.08cm}\V\|_F.\label{cont11}
\end{eqnarray}
\end{lemma}
\begin{proof}
For the first inequality, we have
\begin{eqnarray}
\begin{aligned}\notag
&\|\U\U^T-\V\V^T\|_F\\
\leq& \|\U\U^T-\U\V^T\|_F+\|\U\V^T-\V\V^T\|_F\\
\leq& \|\U\|_2\|\U-\V\|_F+\|\V\|_2\|\U-\V\|_F\\
=& (\|\U\|_2+\|\V\|_2)\|\U-\V\|_F.
\end{aligned}
\end{eqnarray}
For the second one, we have
\begin{eqnarray}
\|\nabla f(\V\V^T)-f(\U^*{\U^*}^T)\|_F&\leq& L\|\V\V^T-\U^*{\U^*}^T\|_F\notag\\
&\leq& L(\|\V\|_2+\|\U^*\|_2)\|\V-\U^*\|_F,\notag
\end{eqnarray}
where we use (\ref{cont10}). For the third one, we have
 \begin{eqnarray}
\begin{aligned}\notag
&\|\nabla f(\U\U^T)\U-\nabla f(\V\V^T)\V\|_F\\
\leq&\|\nabla f(\U\U^T)\U-\nabla f(\V\V^T)\U\|_F+\|\nabla f(\V\V^T)\U-\nabla f(\V\V^T)\V\|_F\\
\leq&\|\U\|_2\|\nabla f(\U\U^T)-\nabla f(\V\V^T)\|_F+\|\nabla f(\V\V^T)\|_2\|\U-\V\|_F\\
\leq&L\|\U\|_2(\|\U\|_2+\|\V\|_2)\|\U-\V\|_F+\|\nabla f(\V\V^T)\|_2\|\U-\V\|_F,
\end{aligned}
\end{eqnarray}
where we use the restricted smoothness of $f$ and (\ref{cont10}) in the last inequality.
\end{proof}
Now we give the proof of Corollary \ref{Lipschitz_constant_lemma}.
\begin{proof}
From Lemma \ref{variable_bound} and the assumptions, we have
\begin{eqnarray}
&\|\U-\U^*\|_F\leq 0.01\sigma_r(\U^*),\notag\\
&\|\U_{ S}-\U_{ S}^*\|_F\leq 0.01\sigma_r(\U_{ S}^*),\notag\\
&0.99\sigma_r(\U^*)\leq \sigma_r(\U)\leq 1.01\sigma_r(\U^*),\notag\\
&0.99\|\U^*\|_2\leq \|\U\|_2\leq 1.01\|\U^*\|_2,\notag\\
&0.99\sigma_r(\U_{ S}^*)\leq \sigma_r(\U_{ S})\leq 1.01\sigma_r(\U_{ S}^*),\notag\\
&0.99\|\U_{ S}^*\|_2\leq \|\U_{ S}\|_2\leq 1.01\|\U_{ S}^*\|_2,\notag
\end{eqnarray}
where $\U$ can be $\U^k$, $\V^k$ and $\Z^k$. From (\ref{cont12}), we have
\begin{eqnarray}
\begin{aligned}
&\|\nabla f(\V^k(\V^k)^T)\|_2\leq \|\nabla f(\V^k(\V^k)^T)-\nabla f(\X^*)\|_2+\|\nabla f(\X^*)\|_2\\
&\leq2.01L\|\U^*\|_2\|\V^k-\U^*\|_F+\|\nabla f(\X^*)\|_2\leq 0.0201L\|\U^*\|_2^2+\|\nabla f(\X^*)\|_2,\label{gradient_bound2}
\end{aligned}
\end{eqnarray}
where we use $\|\V^k-\U^*\|_F\leq 0.01\|\U^*\|_2$. On the other hand, let
\begin{eqnarray}
\hat\Z^{k+1}=\Z^k-\frac{\eta}{\theta_k}\nabla g(\V^k)\notag
\end{eqnarray}
then we have $\Z^{k+1}=\mbox{Project}_{\Omega_{ S}}(\hat\Z^{k+1})$ and
\begin{eqnarray}
\|\hat\Z^{k+1}\|_2&\leq& \|\Z^k\|_2+\frac{2\eta}{\theta_k}\|\nabla f(\V^k(\V^k)^T)\|_2\|\V^k\|_2\notag\\
&\leq&1.01\|\U^*\|_2\left(1+\frac{(0.0402L\|\U^*\|_2^2+2\|\nabla f(\X^*)\|_2)\eta}{\theta_k}\right)\notag\\
&\leq&1.01\|\U^*\|_2\left(1+\frac{1}{\theta_k}\right),\notag
\end{eqnarray}
where we use $\|\Z^k\|_2\leq 1.01\|\U^*\|_2$, $\|\V^k\|_2\leq 1.01\|\U^*\|_2$, (\ref{gradient_bound2}) and the setting of $\eta$. Let $\hat\Omega_{ S}=\{\U_{ S}\in\mathcal{R}^{r\times r}:\U_{ S}\succeq \epsilon \I\}$, then
\begin{eqnarray}
&&\mbox{Project}_{\hat\Omega_{ S}}(\hat\Z_{ S}^{k+1})\notag\\
&=&\argmin_{\U\in\hat\Omega_{ S}}\|\U-\hat\Z_{ S}^{k+1}\|_F^2\notag\\
&=&\argmin_{\U\in\hat\Omega_{ S}}\left\|\U-\frac{\hat\Z_{ S}^{k+1}+(\hat\Z_{ S}^{k+1})^T}{2}-\frac{\hat\Z_{ S}^{k+1}-(\hat\Z_{ S}^{k+1})^T}{2}\right\|_F^2\notag\\
&=&\argmin_{\U\in\hat\Omega_{ S}}\left\|\U-\frac{\hat\Z_{ S}^{k+1}+(\hat\Z_{ S}^{k+1})^T}{2}\right\|_F^2+\left\|\frac{\hat\Z_{ S}^{k+1}-(\hat\Z_{ S}^{k+1})^T}{2}\right\|_F^2,\notag
\end{eqnarray}
where we use $\tr(\A\B)=0$ if $\A=\A^T$ and $\B=-\B^T$, and $\U=\U^T$ from $\U\in\hat\Omega_{ S}$. Let $\U\Sigma\U^T$ be the eigenvalue decomposition of $\frac{\hat\Z_{ S}^{k+1}+(\hat\Z_{ S}^{k+1})^T}{2}$ and $\hat\Sigma_{i,i}=\max\{\epsilon,\Sigma_{i,i}\}$. Then $\mbox{Project}_{\hat\Omega_{ S}}(\hat\Z_{ S}^{k+1})=\U\hat\Sigma\U^T$ and
\begin{eqnarray}
&&\|\mbox{Project}_{\hat\Omega_{ S}}(\hat\Z_{ S}^{k+1})\|_2=\max\{\epsilon,\Sigma_{1,1}\}\notag\\
&&\leq \max\left\{\epsilon,\left\|\frac{\hat\Z_{ S}^{k+1}+(\hat\Z_{ S}^{k+1})^T}{2}\right\|_2\right\}\leq \max\left\{\epsilon,\|\hat\Z_{ S}^{k+1}\|_2\right\},\notag
\end{eqnarray}
where $\Sigma_{1,1}$ is the largest eigenvalue of $\frac{\hat\Z_{ S}^{k+1}+(\hat\Z_{ S}^{k+1})^T}{2}$. Let $\Z_{- S}$ be the submatrix with the rows indicated by the indexes out of $ S$, then
\begin{eqnarray}
\|\Z^{k+1}\|_2&\leq&\|\Z_{ S}^{k+1}\|_2+\|\Z_{- S}^{k+1}\|_2\notag\\
&=&\|\mbox{Project}_{\hat\Omega_{ S}}(\hat\Z_{ S}^{k+1})\|_2+\|\hat\Z_{- S}^{k+1}\|_2\notag\\
&\leq& \max\left\{\epsilon,\|\hat\Z_{ S}^{k+1}\|_2\right\}+\|\hat\Z_{- S}^{k+1}\|_2\notag\\
&\leq& \max\left\{\epsilon,\|\hat\Z^{k+1}\|_2\right\}+\|\hat\Z^{k+1}\|_2\notag\\
&\leq&2\max\left\{\epsilon,\|\hat\Z^{k+1}\|_2\right\}\notag
\end{eqnarray}
and
\begin{eqnarray}
\|\U^{k+1}\|_2&\leq& (1-\theta_{k})\|\U^k\|_2+\theta_{k}\|\Z^{k+1}\|_2\notag\\
&\leq& 1.01(1-\theta_{k})\|\U^*\|_2+2\theta_{k}\max\left\{\epsilon,1.01\|\U^*\|_2\left(1+\frac{1}{\theta_k}\right)\right\}\notag\\
&\leq& 1.01(1-\theta_{k})\|\U^*\|_2+\max\left\{2\epsilon,1.01\|\U^*\|_2\left(2+2\right)\right\}\notag\\
&\leq& 5.05\|\U^*\|_2,\notag
\end{eqnarray}
where we use (\ref{U_step}) in the first inequality, $0\leq\theta_k\leq 1$ in the third and forth inequality and $\|\U^*\|_2\geq\|\U^*_{ S}\|_2\geq \sigma_r(\U^*_{ S})\geq \epsilon$ in the last inequality. So
\begin{eqnarray}
&&\|\nabla f(\V^k(\V^k)^T)\|_2+\frac{L(\|\V^k\|_2+\|\U^{k+1}\|_2)^2}{2}\notag\\
&\leq& 0.0201L\|\U^*\|_2^2+\|\nabla f(\X^*)\|_2+\frac{L(6.06\|\U^*\|_2)^2}{2}\leq\frac{L_g}{2}.\notag
\end{eqnarray}
From Theorem \ref{Lipschitz_S}, we can have the conclusion.
\end{proof}

\section*{Appendix B}\label{sec_appendixB}

\begin{lemma}\label{strong_CW}
Assume that $\U^*\in\mathcal{X}^*$. Then, for any $\U$, we have
\begin{eqnarray}
g(\U)-g(\U^*)\geq 0.4\mu\sigma_r^2(\U^*)\|P_{\mathcal{X}^*}(\U)-\U\|_F^2.\notag
\end{eqnarray}
\end{lemma}
\begin{proof}
From (\ref{cont1}), we have
\begin{eqnarray}
&&f(\U^*{\U^*}^T)-f(\U\U^T)\notag\\
&\leq&2\<\nabla f(\U^*{\U^*}^T)\U^*,\U^*-\U\>-\<\nabla f(\U^*{\U^*}^T),(\U^*-\U)(\U^*-\U)^T\>\notag\\
&&-\frac{\mu}{2}\|\U^*{\U^*}^T-\U\U^T\|_F^2.\notag
\end{eqnarray}
Since $\U^*$ is a minimizer of problem (\ref{non_convex_problem0}), we have $\nabla f(\U^*{\U^*}^T)\U^*=0$. From

\noindent$\<\nabla f(\U^*{\U^*}^T),(\U^*-\U)(\U^*-\U)^T\>\geq0$ and Lemma \ref{diff_lower_bound}, we can have the conclusion.

\end{proof}
Now we give the proof of Theorem \ref{conver_ratem}.
\begin{proof}
Let $\U^{t,*}\in\Omega_{S}\cap\mathcal{X}^*$, where $S=S^1$ when $t$ is odd and $S=S^2$ when $t$ is even. Specially, $\U^{0,*}=\U^*$. From (\ref{aaa3}) we have
\begin{eqnarray}
&&\|\U^{t+1,0}-\U^{t+1,*}\|_F\notag\\
&\leq& \left(\frac{1}{4}\right)^{t+1}\|\U^{0,0}-\U^{0,*}\|_F\notag\\
&=& \left(4^{-\frac{\sqrt{\eta\mu}\sigma_r(\U^{*})\min\{\sigma_r(\U_{ S^1}^{*}),\sigma_r(\U_{ S^2}^{*})\}}{28\|\U^{*}\|_2}}\right)^{(t+1)(K+1)}\|\U^{0,0}-\U^{0,*}\|_F\notag\\
&\leq& \left(1-\frac{\sqrt{\eta\mu}\sigma_r(\U^{*})\min\{\sigma_r(\U_{ S^1}^{*}),\sigma_r(\U_{ S^2}^{*})\}}{28\|\U^*\|_2}\right)^{(t+1)(K+1)}\|\U^{0,0}-\U^{0,*}\|_F,\notag
\end{eqnarray}
where we use $4^{-x}\leq e^{-x}\leq 1-x$.

From Theorem \ref{Lipschitz_S} and $\nabla g(\U^{t+1,*})=0$, we have
\begin{eqnarray}
&&g(\U^{t+1,0})-g(\U^{*})=g(\U^{t+1,0})-g(\U^{t+1,*})\leq\frac{1}{\eta}\|\U^{t+1,0}-\U^{t+1,*}\|_F^2,\notag
\end{eqnarray}
which leads to the conclusion.

\end{proof}

\section*{Appendix C}\label{sec_appendixC}

Proof of Lemma \ref{global_lemma}.
\begin{proof}
We can easily check that $\beta_{\max}<1$ due to $\beta_k\leq 1-\theta_{k-1}$ and the fact that $K$ a finite constant. From Theorem \ref{Lipschitz_S}, we have
\begin{eqnarray}
g(\U^{k+1})&\leq& g(\U^k)+\<\nabla g(\U^k),\U^{k+1}-\U^k\>+\frac{\hat L}{2}\|\U^{k+1}-\U^k\|_F^2\notag\\
&=& g(\U^k)+\<\nabla g(\U^k)-\nabla g(\V^k),\U^{k+1}-\U^k\>\notag\\
&&+\<\nabla g(\V^k),\U^{k+1}-\U^k\>+\frac{\hat L}{2}\|\U^{k+1}-\U^k\|_F^2.\notag
\end{eqnarray}
Applying the inequality of $\<\u,\v\>\leq \|\u\|\|\v\|$, Lemma \ref{matrix_lemma} and the inequality of $2\|\u\|\|\v\|\leq \alpha\|\u\|^2+\frac{1}{\alpha}\|\v\|^2$ to the second term, we can have
\begin{eqnarray}
g(\U^{k+1})&\leq& g(\U^k)+\frac{\hat L}{2}\left(\alpha\|\U^k-\V^k\|_F^2+\frac{1}{\alpha}\|\U^{k+1}-\U^k\|_F^2\right)\notag\\
&&+\<\nabla g(\V^k),\U^{k+1}-\U^k\>+\frac{\hat L}{2}\|\U^{k+1}-\U^k\|_F^2.\notag
\end{eqnarray}
Applying Lemma \ref{temp_lemma} in Appendix C to bound the third term, we can have
\begin{eqnarray}
\begin{aligned}\label{cont5}
&g(\U^{k+1})-g(\U^k)\\
\leq& \frac{\hat L\alpha}{2}\|\U^k-\V^k\|_F^2+\frac{\hat L}{2\alpha}\|\U^{k+1}-\U^k\|_F^2+\frac{1}{2\eta}\|\U^k-\V^k\|_F^2\\
&-\frac{1}{2\eta}\|\U^k-\U^{k+1}\|_F^2+\frac{\hat L}{2}\|\U^{k+1}-\U^k\|_F^2\\
=& \beta_k^2\left(\frac{1}{2\eta}+\frac{\hat L\alpha}{2}\right)\|\U^k-\U^{k-1}\|_F^2-\left(\frac{1}{2\eta}-\frac{\hat L}{2}-\frac{\hat L}{2\alpha}\right)\|\U^{k+1}-\U^k\|_F^2
\end{aligned}
\end{eqnarray}
for all $k=1,2,\cdots,K$, where we use $\V^k-\U^k=\beta_k(\U^k-\U^{k-1})$ proved in Lemma \ref{temp_lemma}. Specially, from {$\U^0=\V^0$} we have
\begin{eqnarray}
g(\U^1)\leq g(\U^0)-\left(\frac{1}{2\eta}-\frac{\hat L}{2}-\frac{\hat L}{2\alpha}\right)\|\U^1-\U^0\|_F^2.\label{cont6}
\end{eqnarray}
Summing (\ref{cont5}) over $k=1,2,\dots,K$ and (\ref{cont6}), we have
\begin{eqnarray}
&&g(\U^{K+1})-g(\U^0)\notag\\
&\leq&-\sum_{k=0}^K\left(\left(\frac{1}{2\eta}-\frac{\hat L}{2}-\frac{\hat L}{2\alpha}\right)-\beta_{k+1}^2\left(\frac{1}{2\eta}+\frac{\hat L\alpha}{2}\right)\right)\|\U^{k+1}-\U^k\|_F^2.\notag
\end{eqnarray}
Letting $\alpha=1/\beta_{\max}$, from the setting of $\eta$, we have the desired conclusion.
\end{proof}

\begin{lemma}\label{temp_lemma}
For Algorithm \ref{agd_alg}, we have
\begin{eqnarray}
\<\nabla g(\V^k),\U^{k+1}-\U^k\>\leq\frac{1}{2\eta}\|\U^k-\V^k\|_F^2-\frac{1}{2\eta}\|\U^k-\U^{k+1}\|_F^2.\notag
\end{eqnarray}
and
\begin{eqnarray}
\V^k=\U^k+\beta_k(\U^k-\U^{k-1}).\notag
\end{eqnarray}
\end{lemma}
\begin{proof}
Let $I_{\Omega_{ S}}(\U)$ be the indicator function of $\Omega_{ S}$. Then, from the optimality condition of (\ref{Z_step}), we have
\begin{eqnarray}
0\in \frac{\theta_k}{\eta}\left(\Z^{k+1}-\Z^k\right)+\nabla g(\V^k)+\partial I_{\Omega_{ S}}(\Z^{k+1}).\notag
\end{eqnarray}
Since $\Omega_{ S}$ is a convex set, we have
\begin{eqnarray}
&&I_{\Omega_{ S}}(\U)\geq I_{\Omega_{ S}}(\Z^{k+1})-\<\frac{\theta_k}{\eta}\left(\Z^{k+1}-\Z^k\right)+\nabla g(\V^k),\U-\Z^{k+1}\>,\forall \U\in\Omega_{ S}\notag
\end{eqnarray}
and
\begin{eqnarray}\label{convex_indicator}
\begin{aligned}
&\frac{\theta_k}{\eta}\<\Z^{k+1}-\Z^k,\U-\Z^{k+1}\>\geq -\<\nabla g(\V^k),\U-\Z^{k+1}\>,\forall \U\in\Omega_{ S}.
\end{aligned}
\end{eqnarray}
With some simple computations, we have
\begin{eqnarray}
&&\<\nabla g(\V^k),\U^{k+1}-\U^k\>\notag\\
&=&\theta_k\<\nabla g(\V^k),\Z^{k+1}-\U^k\>(\mbox{ from } (\ref{U_step}))\notag\\
&\leq&\frac{\theta_k^2}{\eta}\<\Z^{k+1}-\Z^k,\U^k-\Z^{k+1}\>(\mbox{ from } (\ref{convex_indicator}))\notag\\
&=&\frac{\theta_k}{\eta}\<\Z^{k+1}-\Z^k,\U^k-\U^{k+1}\>(\mbox{ from } (\ref{U_step}))\notag\\
&=&\frac{1}{\eta}\<\U^{k+1}-\V^k,\U^k-\U^{k+1}\>(\mbox{ from } (\ref{V_step})\mbox{ and } (\ref{U_step}))\notag\\
&=&\frac{1}{2\eta}\left[\|\U^k-\V^k\|_F^2-\|\U^{k+1}-\V^k\|_F^2-\|\U^k-\U^{k+1}\|_F^2\right]\notag\\
&\leq&\frac{1}{2\eta}\|\U^k-\V^k\|_F^2-\frac{1}{2\eta}\|\U^k-\U^{k+1}\|_F^2.\notag
\end{eqnarray}
From (\ref{V_step}) and (\ref{U_step}), we have
\begin{eqnarray}
\V^k\hspace*{-0.08cm}=\hspace*{-0.08cm}(1\hspace*{-0.08cm}-\hspace*{-0.08cm}\theta_k)\U^k\hspace*{-0.08cm}+\hspace*{-0.08cm}\frac{\theta_k}{\theta_{k-1}}(\U^k\hspace*{-0.08cm}-\hspace*{-0.08cm}(1\hspace*{-0.08cm}-\hspace*{-0.08cm}\theta_{k-1})\U^{k-1})\hspace*{-0.08cm}=\hspace*{-0.08cm}\U^k\hspace*{-0.08cm}+\hspace*{-0.08cm}\frac{\theta_k(1\hspace*{-0.08cm}-\hspace*{-0.08cm}\theta_{k-1})}{\theta_{k-1}}(\U^k\hspace*{-0.08cm}-\hspace*{-0.08cm}\U^{k-1}),\notag
\end{eqnarray}
which leads to the second conclusion.
\end{proof}

\begin{lemma}\label{temp_lemma2}
Under the assumptions in Theorem \ref{global_convergencem}, if (\ref{cont7}) holds, then we have $\|\U^{t',k+1}-\Z^{t',k+1}\|_F\leq\frac{2\varepsilon}{\theta_k}$, $\|\Z^{t',k+1}-\Z^{t',k}\|_F\leq\frac{5\varepsilon}{\theta_k}$ and $\|\Z^{t',k+1}-\V^{t',k}\|_F\leq\frac{9\varepsilon}{\theta_k}$ for $t'=t$ or $t'=t+1$.
\end{lemma}
\begin{proof}
From (\ref{cont7}), for $t'=t$ or $t+1$ and $\forall k=0,\cdots,K$, we can have the following easy-to-check inequalities.
\begin{eqnarray}
&&\U^{t',0}=\V^{t',0}=\Z^{t',0},\label{lim0}\\
&&\|\U^{t',k+1}-\U^{t',k}\|_F\leq\varepsilon,\label{lim1}\\
&&\|\Z^{t',k+1}-\U^{t',k}\|_F\leq\frac{\varepsilon}{\theta_k},(\mbox{ from } (\ref{U_step}))\label{lim2}\\
&&\|\U^{t',k+1}-\Z^{t',k+1}\|_F\leq\varepsilon+\frac{\varepsilon}{\theta_k},(\mbox{ from } (\ref{lim1})\mbox{ and }(\ref{lim2}))\label{lim3}\\
&&\|\V^{t',k+1}-\U^{t',k+1}\|_F\leq\theta_{k+1}\left(\varepsilon+\frac{\varepsilon}{\theta_{k}}\right),(\mbox{ from } (\ref{V_step})\mbox{ and } (\ref{lim3}))\label{lim4}\\
&&\|\Z^{t',k+1}-\Z^{t',k}\|_F\leq \varepsilon+\frac{\varepsilon}{\theta_k}+\varepsilon+\frac{\varepsilon}{\theta_{k-1}}+\varepsilon,(\mbox{ from } (\ref{lim1}) \mbox{ and } (\ref{lim3}))\label{lim5}\\
&&\|\Z^{t',k+1}-\V^{t',k}\|_F\leq \varepsilon\hspace*{-0.08cm}+\hspace*{-0.08cm}\frac{\varepsilon}{\theta_k}\hspace*{-0.08cm}+\hspace*{-0.08cm}(2\hspace*{-0.08cm}+\hspace*{-0.08cm}\theta_k)\left(\varepsilon\hspace*{-0.08cm}+\hspace*{-0.08cm}\frac{\varepsilon}{\theta_{k-1}}\right)\hspace*{-0.08cm}+\hspace*{-0.08cm}\varepsilon,( (\ref{lim5}),(\ref{lim3}),(\ref{lim4}))\label{lim6}
\end{eqnarray}
From $\theta_k\leq\theta_{k-1}\leq 1$, we can have the conclusions.
\end{proof}

\section*{Appendix D}\label{appendixD}

\begin{lemma}\label{strong_C2}
Assume that $\U^*\in\mathcal{X}^*$ and $\V\in\mathcal{R}^{n\times r}$ satisfy $\|\V-P_{\mathcal{X}^*}(\V)\|_F\leq \min\left\{0.01\sigma_r(\U^*), \frac{\mu\sigma_r^2(\U^*)}{6L\|\U^*\|_2}\right\}$. Then, we have
\begin{eqnarray}
\|\V-P_{\mathcal{X}^*}(\V)\|_F\leq\frac{5}{\mu\sigma_r^2(\U^*)}\|\nabla g(\V)\|_F.\notag
\end{eqnarray}
\end{lemma}
\begin{proof}
Similar to the proof of Theorem \ref{strong_Cm}, we have
\begin{eqnarray}
\begin{aligned}\label{cont13}
&g(\U^*)=g(P_{\mathcal{X}^*}(\V))\\
&\geq g(\V)+\<\nabla g(\V),P_{\mathcal{X}^*}(\V)-\V\>+0.2\mu\sigma_r^2(\U^*)\|\V-P_{\mathcal{X}^*}(\V)\|_F^2,
\end{aligned}
\end{eqnarray}
where we use Lemma \ref{diff_lower_bound} to bound $\|\V\V^T-P_{\mathcal{X}^*}(\V)(P_{\mathcal{X}^*}(\V))^T\|_F^2$. Since $g(\U^*)\leq g(\V)$, we can have
\begin{eqnarray}
0.2\mu\sigma_r^2(\U^*)\|\V-P_{\mathcal{X}^*}(\V)\|_F^2&\leq& \<\nabla g(\V),\V-P_{\mathcal{X}^*}(\V)\>\notag\\
&\leq&\|\nabla g(\V)\|_F\|\V-P_{\mathcal{X}^*}(\V)\|_F,\notag
\end{eqnarray}
which leads to the conclusion.
\end{proof}

\begin{lemma}\label{decrease_lemma}
Under the assumptions of Lemma \ref{global_lemma}, we have
\begin{eqnarray}
&&g(\U^{k+1})+\nu\|\U^{k+1}-\U^k\|_F^2\notag\\
&\leq& g(\U^k)+\nu\|\U^k-\U^{k-1}\|_F^2-\gamma\left(\|\U^{k}-\U^{k-1}\|_F^2+\|\U^{k+1}-\U^k\|_F^2\right),\notag
\end{eqnarray}
where $\gamma=\frac{1-\beta_{\max}^2}{4\eta}-\frac{\beta_{\max}\hat L}{2}-\frac{\hat L}{4}>0$ and $\nu=\frac{1+\beta_{\max}^2}{4\eta}-\frac{\hat L}{4}>0$.
\end{lemma}
\begin{proof}
Letting $\alpha=\frac{1}{\beta_{\max}}$ in (\ref{cont5}), we can have the conclusion.
\end{proof}
Now we give the proof of Theorem \ref{direct_min_theorem}.
\begin{proof}
Denote $\hat\U^*=P_{\mathcal{X}^*}(\V^k)$. From Theorem \ref{Lipschitz_S}, we can have
\begin{eqnarray}
&&g(\U^{k+1})\notag\\
&\leq& g(\V^k)+\<\nabla g(\V^k),\U^{k+1}-\V^k\>+\frac{\hat L}{2}\|\U^{k+1}-\V^k\|_F^2\notag\\
&=& g(\V^k)+\<\nabla g(\V^k),\hat\U^*-\V^k\>+\<\nabla g(\V^k),\U^{k+1}-\hat\U^*\>+\frac{\hat L}{2}\|\U^{k+1}-\V^k\|_F^2\notag\\
&\leq& g(\hat\U^*)+\<\nabla g(\V^k),\U^{k+1}-\hat\U^*\>+\frac{\hat L}{2}\|\U^{k+1}-\V^k\|_F^2\notag\\
&=& g(\hat\U^*)+\frac{1}{\eta}\<\V^k-\U^{k+1},\U^{k+1}-\hat\U^*\>+\frac{\hat L}{2}\|\U^{k+1}-\V^k\|_F^2\notag\\
&\leq& g(\hat\U^*)+\frac{1}{\eta}\<\V^k-\U^{k+1},\V^k-\hat\U^*\>\notag\\
&\leq& g(\hat\U^*)+\frac{1}{\eta}\|\V^k-\U^{k+1}\|_F\|\V^k-\hat\U^*\|_F\notag\\
&\leq& g(\hat\U^*)+\frac{5}{\eta\mu\sigma_r^2(\U^*)}\|\V^k-\U^{k+1}\|_F\|\nabla g(\V^k)\|_F\notag\\
&=& g(\hat\U^*)+\frac{5}{\eta^2\mu\sigma_r^2(\U^*)}\|\U^{k+1}-\V^k\|_F^2\notag\\
&=& g(\hat\U^*)+\frac{5}{\eta^2\mu\sigma_r^2(\U^*)}\|\U^{k+1}-\U^k-\beta_k(\U^k-\U^{k-1})\|_F^2\notag\\
&\leq& g(\hat\U^*)+\frac{5}{\eta^2\mu\sigma_r^2(\U^*)}\left(\|\U^{k+1}-\U^k\|_F^2+\|\U^k-\U^{k-1}\|_F^2\right),\notag
\end{eqnarray}
where we use (\ref{cont13}) in the second inequality, (\ref{U_step3}) in the second equality,  $\eta<\frac{1}{\hat L}$ in the third inequality, Lemma \ref{strong_C2} in the forth inequality, (\ref{V_step3}) in the fifth equality and $\beta_{\max}<1$ in the last inequality. So we have
\begin{eqnarray}
\begin{aligned}\label{cont14}
&g(\U^{k+1})+\nu\|\U^{k+1}-\U^k\|_F^2-g(\U^*)\\
\leq&\left(\frac{5}{\eta^2\mu\sigma_r^2(\U^*)}+\nu\right)\left(\|\U^{k+1}-\U^k\|_F^2+\|\U^k-\U^{k-1}\|_F^2\right).
\end{aligned}
\end{eqnarray}
Combing Lemma \ref{decrease_lemma} and (\ref{cont14}), we can have the conclusion.
\end{proof}

\small
\bibliographystyle{unsrt}
\bibliography{agd}

\end{document}